\documentclass[preprint,authoryear,11pt,a4paper]{elsarticle}

\usepackage{geometry}
\usepackage{setspace}
\usepackage{apalike}

\usepackage[hyphens]{url}
\usepackage[caption=false]{subfig}
\usepackage{float}
\usepackage{hyperref}
\usepackage{algorithm, algorithmic}
\usepackage{amsmath,amssymb,amsfonts,amsthm}
\allowdisplaybreaks[4]
\usepackage{graphicx}
\usepackage{enumerate}
\usepackage{array}
\usepackage{multirow}
\usepackage{color}

\usepackage{tikz}

\newtheorem{theorem}{Theorem}
\newtheorem{lemma}{Lemma}
\newtheorem{corollary}{Corollary}
\newtheorem{proposition}{Proposition}

\newtheorem{definition}{Definition}

\newtheorem{example}{Example}

\journal{European Journal of Operational Research}

\begin{document}
\onehalfspacing

\begin{frontmatter}
\title{Multistage Robust Optimization for Time-Decoupled Power Flexibility Aggregation with Energy Storage}

\author[Address1]{Rui Xie}
\ead{ruixie@cuhk.edu.hk}
\author[Address1]{Yue Chen\corref{cor1}}
\ead{yuechen@mae.cuhk.edu.hk}
\cortext[cor1]{Corresponding author}
\address[Address1]{Department of Mechanical and Automation Engineering, The Chinese University of Hong Kong, Hong Kong.}

\begin{abstract} 
To mitigate global climate change, distributed energy resources (DERs), such as distributed generators, flexible loads, and energy storage systems (ESSs), have witnessed rapid growth in power distribution systems. When properly managed, these DERs can provide significant flexibility to power systems, enhancing both reliability and economic efficiency. Due to their relatively small scale, DERs are typically managed by the distribution system operator (DSO), who interacts with the transmission system operator (TSO) on their behalf. Specifically, the DSO aggregates the power flexibility of the DERs under its control, representing it as a feasible variation range of aggregate active power at the substation level. This flexibility range is submitted to the TSO, who determines a setpoint within that range. The DSO then disaggregates the setpoint to dispatch DERs. This paper focuses on the DSO’s power flexibility aggregation problem. First, we propose a novel multistage robust optimization model with decision-dependent uncertainty for power flexibility aggregation. Distinct from the traditional two-stage models, our multistage framework captures the sequential decision-making of the TSO and DSO and is more general (e.g., can accommodate non-ideal ESSs). Then, we develop multiple solution methods, including exact, inner, and outer approximation approaches under different assumptions, and compare their performance in terms of applicability, optimality, and computational efficiency. Furthermore, we design greedy algorithms for DSO's real-time disaggregation. We prove that the rectangular method yields greater total aggregate flexibility compared to the existing approach. Case studies demonstrate the effectiveness of the proposed aggregation and disaggregation methods, validating their practical applicability.
\end{abstract}

\begin{keyword}
OR in energy, multistage robust optimization, energy storage, power flexibility aggregation, decision-dependent uncertainty
\end{keyword}

\end{frontmatter}

\section{Introduction}

To mitigate global warming, numerous countries and regions have established plans to achieve carbon neutrality. The energy sector plays a pivotal role in this transition, particularly through renewable energy integration and energy efficiency improvements. In this context, distributed energy resources (DERs), such as distributed generators (DGs), flexible loads, and energy storage systems (ESSs), have grown rapidly in power distribution systems over the past few decades \citep{gust2024designing}. This growth has introduced substantial flexibility at the distribution power systems to support the upstream transmission system's energy management and enhance the efficiency and reliability of the overall power system. However, DERs cannot participate directly in the transmission-level operation optimization due to their small capacities. Fortunately, distribution system operators (DSOs) can aggregate the distributed flexibility and interact with the transmission system operator (TSO) on behalf of the DERs. The DER aggregators in the U.S. have been permitted to participate in transmission-level energy markets since 2020 by the FERC Order 2222
\citep{federalregister2021}. 

Due to concerns about computational complexity and privacy leakage, the DSO does not share individual models of DERs with the TSO. Instead, it calculates and submits the feasible variation region of aggregate active power to the TSO as a representation of the total flexibility, which is a process called \emph{aggregation}. Then, the TSO can efficiently utilize this region to form a constraint characterizing the distribution system's aggregate power in the transmission-level optimization. Once the TSO determines a setpoint of the aggregate power for the distribution system, the DSO coordinates individual DERs to realize this setpoint, a process referred to as \emph{disaggregation}. This paper studies methods for DSOs to conduct aggregation and disaggregation processes, with the goal of maximizing aggregate flexibility while minimizing operation costs.


Power flexibility aggregation can be viewed as a dimension reduction process, transforming the operational constraints of DERs and distribution networks into a simpler set of constraints that characterize the feasible variation region of aggregate active power at the substation that connects the distribution power system to the transmission power system. Power flexibility aggregation has been widely applied in various scenarios, including residential and non-residential load scheduling \citep{ayon2017optimal}, thermostatically controlled load (TCL) scheduling and control \citep{paridari2020flexibility}, electric vehicle (EV) charging station coordination \citep{yan2023real}, and energy market participation \citep{wang2020aggregation}.


When devices with time-coupled operational constraints such as ESSs are present, multiple time periods must be considered jointly in power flexibility aggregation. This is because operational decisions in earlier periods can affect the flexibility available in later periods through the state of the devices, such as the state-of-charge (SoC) of ESSs. In such cases, the aggregate flexibility is represented by a feasible variation region of trajectories over time, which is a subset of the $T$-dimensional Euclidean space with $T$ the number of periods. This region comprises trajectories of aggregate power that can be achieved by the DSO through dispatching DERs. In previous works' settings, the DSO performs the disaggregation process with a known trajectory of aggregate power, where the definition of the aggregate power flexibility region is as follows:


\begin{definition}[Aggregate power flexibility region (two-stage)]
\label{def:two-stage-region}
    A set $S_A \subseteq \mathbb{R}^T$ is the aggregate power flexibility region of DERs if and only if for any aggregate power trajectory $p^A = (p_1^A, \dots, p_T^A) \in S_A$, there exist operation strategies $\alpha_1(p^A), \alpha_2(p^A), \dots, \alpha_T(p^A)$ for DERs so that their aggregate power in any period $t$ is $p_t^A$. Each $\alpha_t(p^A)$ represents the operation strategies in time period $t$, which may depend on the aggregate power trajectory $p^A$ over all periods.
\end{definition}

In Definition~\ref{def:two-stage-region}, the operation strategies may depend on future aggregate power. However, as renewable energy penetration increases in power systems, it becomes increasingly challenging for the TSO to determine the aggregate power trajectories of distribution systems multiple periods in advance. In practice, real-time dispatch is critical to maintain instantaneous power balance \citep{yildiran2023robust}, especially given the variability and uncertainty introduced by renewable energy sources. Consequently, a new definition of aggregate power flexibility is needed to account for the uncertainty of TSO's future dispatch of DERs. To address this issue, this paper proposes a multistage model that captures the sequential decision-making process of the TSO, where the new definition of the aggregate power flexibility region is as follows:

\begin{definition}[Aggregate power flexibility region (multistage)]
\label{def:multistage-region}
    A set $S_A \subseteq \mathbb{R}^T$ is the aggregate power flexibility region of DERs if and only if for any aggregate power trajectory $p^A = (p_1^A, \dots, p_T^A) \in S_A$, there exist operation strategies $\alpha_1(p_{\leq 1}^A), \alpha_2(p_{\leq 2}^A), \dots, \alpha_T(p_{\leq T}^A)$ for the DERs so that their aggregate power in any period $t$ is $p^A$. Each $\alpha_t(p_{\leq t}^A)$ represents the operation strategies in period $t$, which may depend on $p_{\leq t}^A = (p_1^A, p_2^A, \dots, p_t^A)$, the aggregate power trajectory up to period $t$.
\end{definition}

\begin{figure}[!t]
\centering
\includegraphics[width=9.0cm]{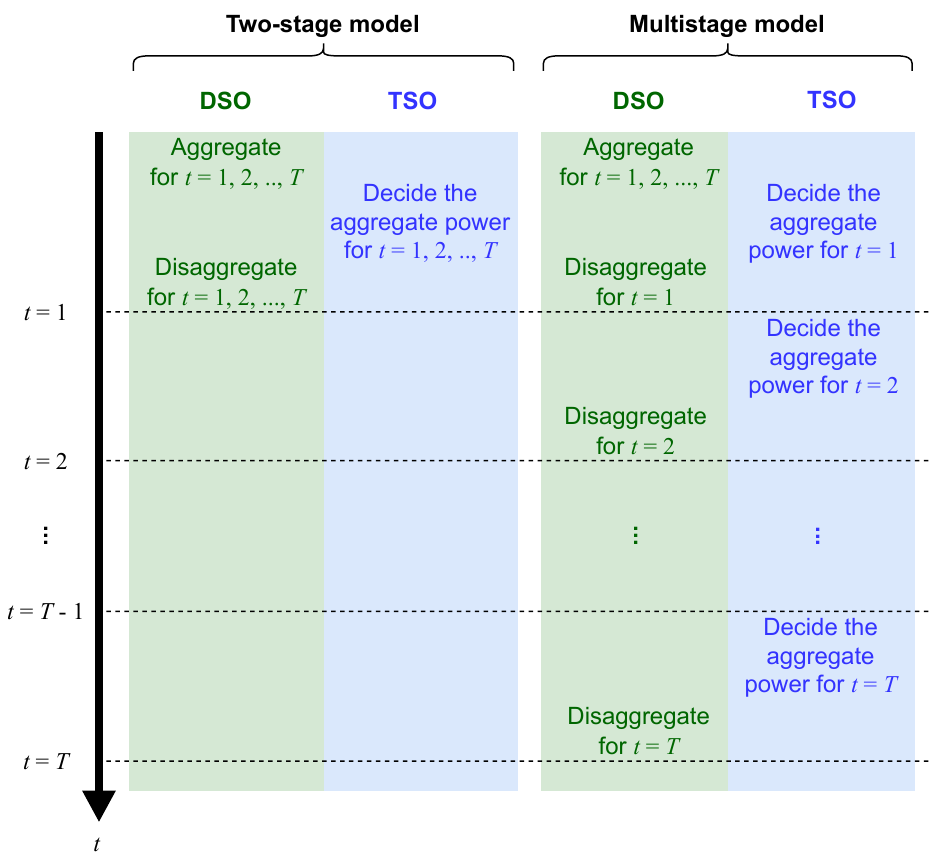}
\caption{Illustration for two-stage and multistage models for aggregate flexibility.}
\label{fig:illustration}
\end{figure}

We compare the proposed multistage model with the traditional approach in Figure~\ref{fig:illustration}. In the traditional framework, the DSO performs disaggregation before period $t=1$ based on the full trajectory of aggregate power over the next $T$ periods, resulting in a two-stage model \citep{chen2021leveraging}. This approach allows the DSO to strategically adjust the disaggregation scheme in earlier periods by leveraging the precise knowledge of future aggregate power trajectories. However, as uncertainty in the transmission system grows, it becomes increasingly challenging for the TSO to accurately determine the entire aggregate power trajectory at the beginning. In period $t \geq 2$, even if the trajectory changes within the power flexibility region, the disaggregation may be infeasible, because the operation strategies before period $t$ cannot be adjusted. This makes the two-stage model impractical. On the contrary, in the proposed multistage framework, the TSO decides the aggregate power and the DSO disaggregates for period $t$ just before period $t$ begins. Thus, the TSO retains more time to decide on aggregate power in response to uncertainties arising in the transmission system.

The two-stage and multistage models described above are fundamentally distinct. Intuitively, in the two-stage model, DSOs can leverage precise knowledge of future aggregate power trajectories for disaggregation, whereas this is not possible in the multistage model. Consequently, the aggregate power in the multistage case must adhere to stricter requirements, resulting in a smaller aggregate power flexibility region compared to the two-stage model. Such additional requirement in the multistage model belongs to nonanticipativity constraint \citep{birge2011introduction}, which ensures that decisions at each stage are made without relying on future realizations of uncertainty, i.e., the future aggregate power trajectory. 

This paper focuses specifically on \emph{time-decoupled} power flexibility aggregation, where the aggregate power flexibility region takes the form of a Cartesian product of $T$ intervals. Formally, we consider the time-decoupled region $S_A = [p_1^{A \vee}, p_1^{A \wedge}] \times [p_2^{A \vee}, p_2^{A \wedge}] \times \dots \times [p_T^{A \vee}, p_T^{A \wedge}]$. An aggregate power trajectory $p^A = (p_1^A, p_2^A, \dots, p_T^A) \in S_A$ if and only if the components satisfy $p_t^{A \vee} \leq p_t^A \leq p_t^{A \wedge}$ for $t = 1, 2, \dots, T$. 
We focus on time-decoupled power flexibility aggregation for the following reasons: 1) The time-decoupled structure is highly beneficial for the TSO in operation management. It enables the TSO to model distribution systems using time-decoupled constraints, which are significantly less complex than time-coupled constraints, such as those involving SoC dynamics. Consequently, time-decoupled power flexibility aggregation has been widely adopted in various applications in the two-stage framework \citep{yan2023real,li2025aggregate}. 2) The time-decoupled region aligns naturally with our multistage framework, represented as the Cartesian product of intervals $[p_t^{A \vee}, p_t^{A \wedge}]$ for $t = 1, 2, \dots, T$, where each interval corresponds to the range of the aggregate power in a period. 3) The time-decoupled structure significantly simplifies the considered model. This simplification will be leveraged later to derive efficient solution methods.

There exist numerous time-decoupled aggregate power flexibility regions for the same group of DERs, which do not necessarily exhibit inclusion relationships with each other. To exploit flexibility, we maximize a flexibility index $\phi(p^{A \vee}, p^{A \wedge})$ over all possible time-decoupled aggregate power flexibility regions, where $\phi: \mathbb{R}^T \times \mathbb{R}^T \rightarrow \mathbb{R}$ is a function. From the DSO's perspective, the operation strategy of the TSO is a kind of uncertainty. Thus, the DSO's flexibility maximization problem can be viewed as a multistage robust optimization (RO) problem. The uncertainty is the aggregate power selected by the TSO, whose range is restricted by the DSO's aggregation results. Thus, the uncertainty in the RO model depends on earlier decisions, known as decision-dependent uncertainty (DDU) or endogenous uncertainty \citep{nohadani2018optimization}. 

The concept of time-decoupled power flexibility aggregation was introduced and explored by \cite{chen2019aggregate} and \cite{chen2021leveraging}, but only within the two-stage framework as an inner approximation for aggregate flexibility. To the best of our knowledge, this paper is the first to investigate time-decoupled power flexibility aggregation from the perspective of a multistage framework. While two-stage RO with DDU can be solved in practice using variations of column-and-constraint generation (C\&CG) algorithm \citep{zeng2013solving,zeng2022two}, obtaining optimal solutions for multistage RO is generally much more challenging due to the complexity introduced by the additional nonanticipativity constraints.

To address the aforementioned research gaps, we first establish a multistage RO model for time-decoupled power flexibility aggregation in distribution systems. Next, we design multiple methods to solve the proposed model either exactly or approximately under various assumptions. Corresponding disaggregation methods are also developed to assist DSOs in achieving the aggregate power specified by the TSO while reducing operation costs. Finally, the proposed methods are compared and validated through case studies.

The contributions of this paper are summarized as follows:

1) \emph{Modeling:} We propose a multistage RO model with DDU for time-decoupled power flexibility aggregation in distribution systems. The model incorporates operational constraints of loads, DGs, ESSs, and power flow model. Unlike previous studies based on two-stage models \citep{chen2019aggregate,chen2021leveraging}, which cannot capture the TSO's sequential decision-making process, our approach employs a multistage framework designed to address this limitation. Additionally, the proposed model is more general, as it accommodates non-ideal ESSs (whose charging and discharging efficiencies can be strictly less than $100\%$), any type of convex power flow models, and general flexibility indices (i.e., the form of function $\phi$ is not restricted).

2) \emph{Methods for aggregation:} We develop an exact enumeration-based solution method for aggregation under convex ESS models. For the general case, we prove that the envelope-based method introduced by \cite{chen2019aggregate} adheres to the nonanticipativity constraints in the multistage model, thereby providing an inner approximation solution, even though it was originally designed for the two-stage model. Additionally, we propose a novel inner approximation method using SoC ranges, which is proven to yield larger flexibility compared to the envelope-based method (by up to 29.9\% in the case studies). Furthermore, we introduce an outer approximation method derived from the two-stage model, which also serves as an efficient outer approximation for the multistage model. These methods collectively provide a framework for flexibility aggregation, balancing accuracy and scalability in practical applications.

3) \emph{Methods for disaggregation:} We propose greedy disaggregation algorithms for the aforementioned exact and inner approximate aggregation methods. Unlike previous works \citep{chen2019aggregate}, which do not account for operation costs, our approaches explicitly minimize the DSOs' operation costs in the current period. Case studies demonstrate that this leads to a reduction in the average operation cost by up to 40.3\%.

4) \emph{Impact of ESS complementarity constraint:} The ESS complementarity constraint, which prevents ESSs from charging and discharging simultaneously, is inherently nonconvex and hence often omitted in power system operation models to maintain convexity \citep{xie2023sizing}. On the necessity of considering the ESS complementarity constraint in power flexibility aggregation, we show that disregarding it for non-ideal ESSs can lead to an overestimation of aggregate flexibility, but this constraint is redundant for ideal ESSs. 
For general cases, a convex approximation can be employed to replace the complementarity constraint \citep{shen2020modeling}.

The rest of this paper is organized as follows. The related literature is reviewed in Section~\ref{sec:review}. The proposed multistage RO model is developed in Section~\ref{sec:model}. The aggregation and disaggregation solution methods are introduced and compared in Section~\ref{sec:solution}. The example cases and numerical experiments are presented in Section~\ref{sec:case}. Finally, Section~\ref{sec:conclusion} concludes the paper.

\section{Related Literature}
\label{sec:review}

In this section, we review recent literature relevant to this work, focusing on power flexibility aggregation and multistage RO with DDU.

\subsection{Power Flexibility Aggregation}

Power flexibility aggregation has been extensively studied in the past decade, with the goal of replacing complex operational constraints of massive devices with simplified aggregation models. One way to do this is to establish a virtual device model. For example, 
\cite{hao2014aggregate} modeled the aggregate flexibility of TCLs as a stochastic non-ideal ESS, with capacity and power bounds dependent on TCL parameters, ambient temperature, and target values. Another approximate virtual ESS model was proposed in \cite{zhao2017geometric} using linear programming (LP). The aggregate flexibility of virtual power plants consisting of DERs was characterized by a hybrid generator-ESS model in \cite{wang2021aggregate} to reduce conservativeness compared to standalone ESS approximations.

Another popular approach to power flexibility aggregation focuses on directly characterizing the feasible variation region of aggregate power trajectories. In this paradigm, \cite{churkin2023tracing} modeled the exact aggregate flexibility of DGs using second-order cone programming (SOCP) relaxations of AC optimal power flow in radial distribution networks. However, this framework excludes devices with time-coupled operational constraints such as ESSs, whose exact flexibility aggregation requires computing Minkowski sums of polyhedra \citep{zhao2017geometric}, proven to be NP-hard \citep{wen2022aggregate}.

To address these computational challenges of characterizing the aggregate power flexibility region, researchers have developed approximation methods categorized as either outer approximations (providing necessary conditions) or inner approximations (ensuring sufficient conditions). Several outer approximate models were proposed in \cite{wen2022aggregate} to avoid the exponential complexity of exact aggregate flexibility derived by Fourier-Motzkin elimination. This work was generalized in \cite{wen2022aggregatetemporally} to further consider distribution system security constraints under a linearized power flow model. However, outer approximations may overestimate the flexibility and include impossible aggregate power trajectories.

In contrast, inner approximations ensure feasibility first at the cost of potential conservatism. Mathematical optimization is often involved in finding the best inner approximations. Among inner approximation techniques, \cite{muller2017aggregation} leveraged zonotopes to model DERs' aggregate flexibility by exploiting the computational efficiency of zonotopic Minkowski summation. An LP problem was formulated to find the inner zonotope with the best approximation quality. In \cite{cui2021network}, two-stage RO was formulated to find the maximum-volume ellipsoidal inner approximations. For EV fleets, \cite{al2024efficient} developed LP-driven maximum-volume polyhedral inner approximations.

Time-decoupled power flexibility aggregation offers a structured inner approximation defined as the Cartesian product of feasible variation region in each periods, i.e., $S_A = [p_1^{A \vee}, p_1^{A \wedge}] \times [p_2^{A \vee}, p_2^{A \wedge}] \times \dots \times [p_T^{A \vee}, p_T^{A \wedge}]$. This structure leads to time-decoupled constraints for aggregate power, which is convenient for the TSO to handle in transmission-level operation management. The maximum time-decoupled flexibility aggregation was formulated as a two-stage RO problem with DDU in \cite{chen2021leveraging}, exactly solved by the C\&CG algorithm. A single-stage inner approximation was proposed in \cite{chen2019aggregate} using the envelopes of ideal ESSs' power trajectories, which is much more efficient than the C\&CG algorithm. This paradigm was extended to various application scenarios: An online EV aggregation method was proposed in \cite{yan2023real} based on time-decoupled aggregate flexibility to deal with the uncertain EV arrivals and charging demands. The aggregate flexibility of multi-energy systems was quantified by time-decoupled sets in \cite{li2025aggregate}, considering power, heating, and gas systems.

While the aforementioned approaches rely on physics-based models, data-driven paradigms have been emerging in the field of flexibility aggregation. For example, 
\cite{taheri2022data} established a hybrid model-informed data-driven methodology to characterize time-coupled load flexibility through inner polyhedral representations. \cite{zhang2024learning} proposed a model-free online aggregation technique for EVs using deep reinforcement learning to adaptively capture real-time flexibility.

Despite the aforementioned progress, there are still research gaps in power flexibility aggregation: 
1) In existing methods, disaggregation occurs only after the entire aggregate power trajectory is determined, which fails to accommodate TSO's sequential decisions in uncertain environments. Thus, we need to formulate a multistage model and develop corresponding solution methods for aggregation and disaggregation.
2) Existing time-decoupled aggregation techniques are usually designed for ideal ESSs and certain types of power flow models and flexibility indices. Generalizations to non-ideal ESSs, arbitrary types of convex power flow models, and general flexibility indices will broaden the applicability.

\subsection{Multistage RO with DDU}

Adjustable RO, including two-stage and multistage RO, is generally NP-hard even under decision-independent uncertainty (DIU), continuous decisions, and linear constraints and objectives \citep{guslitser2002uncertainty}. While algorithms such as Benders decomposition and the C\&CG method \citep{zeng2013solving} effectively solve two-stage RO with DIU, they are unsuitable for multistage settings. Consequently, multistage RO often requires approximation techniques. Early work by \cite{ben2004adjustable} introduced affine policies for adjustable RO.
Enhanced affine/quadratic decision rules via copositive programming were studied in \cite{xu2023improved}. Alternative approaches include finite adaptability, where recourse decisions are restricted to a predefined set of policies \citep{bertsimas2010finite}, or the uncertainty set is partitioned iteratively to make decisions based on realized uncertainties \citep{postek2016multistage}. 
In contrast to approximations, \cite{georghiou2019robust} proposed a robust dual dynamic programming algorithm for multistage RO with DIU that converges to exact solutions under specific conditions.

Multistage RO with DIU has diverse applications. For example, \cite{yildiran2023robust} developed a multistage RO framework for economic dispatch under renewable uncertainty, where exact solution can be found with ideal ESSs. \cite{kim2024sample} combined multistage stochastic programming with RO for virtual power plant bidding. \cite{portoleau2024robust} employed robust decision trees for multistage project scheduling. These applications demonstrate the critical role of multistage RO in managing real-world sequential uncertainties.

When DDU presents, prior decisions influence future uncertainty realizations. This increases the coupling between decisions and uncertainties, leading to a higher complexity of multistage RO compared to the DIU case. Research addressing multistage RO with DDU remains limited. \cite{lappas2016multi} designed decision-dependent uncertainty sets for process scheduling applications, while \cite{zhang2020unified} unified DDU modeling frameworks and derived decision rule approximations. Multistage robust mixed-integer optimization with DDU was studied by \cite{feng2021multistage} and solutions were developed using nonlinear and piecewise linear decision rules. 
Multistage RO with DDU was applied to energy dispatch with demand response by \cite{su2022multi}, where affine policies and scenario mapping are used for solution, but non-ideal ESSs were not considered.

In addition, existing general methods for multistage RO with DDU often suffer from over-conservativeness or computational inefficiency. To address these limitations, novel approaches that exploit the specific model characteristics are needed to effectively solve the multistage RO problem with DDU arising from time-decoupled power flexibility aggregation.

\section{Time-Decoupled Power Flexibility Aggregation Model}
\label{sec:model}

This section develops a time-decoupled power flexibility aggregation framework using multistage RO with DDU. We first model the individual components (i.e., load, DG, ESS, and power network) and subsequently formulate the multistage RO problem.

\subsection{Component Models}
\label{sec:model-component}

\begin{figure}[!t]
\centering
\includegraphics[width=8.5cm]{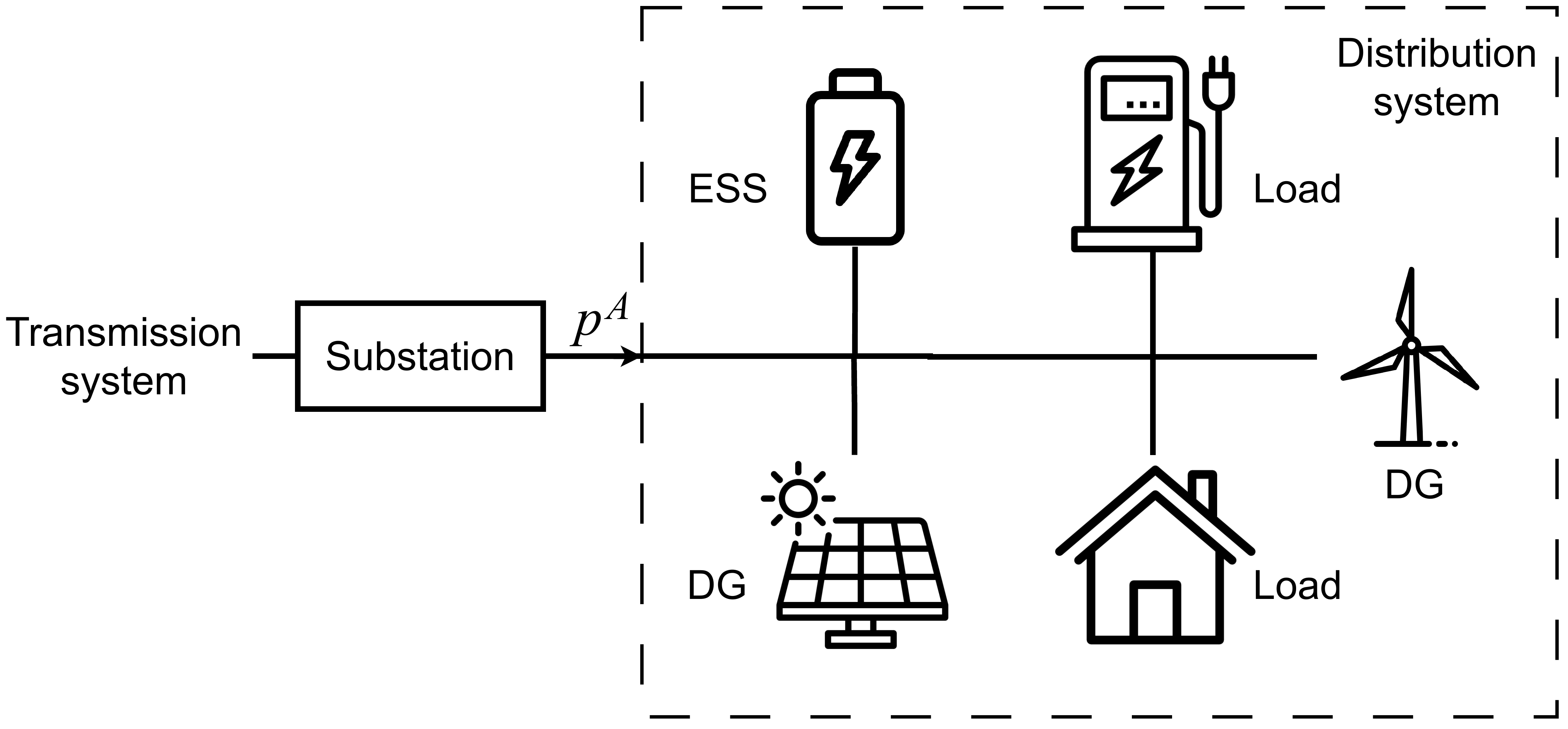}
\caption{Structure of the considered distribution network.}
\label{fig:system}
\end{figure}

We study power flexibility aggregation in a distribution system composed of loads, DGs, ESSs, and transmission lines. The system accommodates both controllable and non-controllable loads or DGs, with its overall structure depicted in Figure~\ref{fig:system}. We consider a finite optimization horizon divided into $T$ periods, indexed by $S_T = \{1, 2, \dots, T\}$. The length of each period is $\tau > 0$. The index set of nodes in the distribution system is denoted by $S_N$, where node $1$ is the reference node and connects the upstream transmission system. We proceed to develop detailed mathematical models for each component.

\subsubsection{Load}

Let $S_D \subseteq S_N$ denote the set of loads. For each load $i \in S_D$ and period $t \in S_T$, let $p_{i t}^D$ and $q_{i t}^D$ represent the active and reactive load power, respectively. The flexibility of these loads is characterized by the following constraints: For any $i \in S_D$ and $t \in S_T$,
\begin{align}
\label{eq:load}
    \underline{P}_{i t}^D \leq p_{i t}^D \leq \overline{P}_{i t}^D,\, q_{i t}^D = \eta_{i t}^D p_{i t}^D.
\end{align}
In the first constraint of \eqref{eq:load}, $\underline{P}_{i t}^D$ and $\overline{P}_{i t}^D$ denote the minimum and maximum active power bounds of load $i$ in period $t$. If load $i$ is uncontrollable, then $\underline{P}_{i t}^D = \overline{P}_{i t}^D$, for any $t \in S_T$. In the second constraint, $\eta_{i t}^D$ is a constant determined by the power factor of load $i$ in period $t$. 

\subsubsection{Distributed Generator}

Denote the set of DGs by $S_G \subseteq S_N$. The active and reactive generation power of DG $i \in S_G$ in period $t \in S_T$ are denoted by $p_{i t}^G$ and $q_{i t}^G$, respectively. The constraints for DGs are formulated as: For any $i \in S_G$ and $t \in S_T$,
\begin{align}
    \label{eq:DER}
    \underline{P}_{i t}^G \leq p_{i t}^G \leq \overline{P}_{i t}^G,\, (p_{i t}^G)^2 + (q_{i t}^G)^2 \leq (\overline{S}_{i t}^G)^2.
\end{align}
The first constraint defines the upper and lower bounds for $p_{i t}^G$, while the second one specifies the maximum apparent power $\overline{S}_{i t}^G$.

\subsubsection{Energy Storage System}

The set of ESSs is denoted by $S_S \subseteq S_N$. For each ESS $i \in S_S$ and period $t \in S_T$, $p_{i t}^S$ and $q_{i t}^S$ represent the active and reactive power injected into the distribution system, respectively. ESS $i$ is discharging when $p_{i t}^S > 0$ and charging when $p_{i t}^S < 0$. The stored energy of ESS $i$ at the end of period $t \in S_T$ is denoted by $e_{i t}^S$, with $e_{i 0}^S$ representing the initial stored energy at the start of the optimization horizon. The constraints for ESSs are as follows: For any $i \in S_S$ and $t \in S_T$,
\begin{subequations}
\label{eq:ESS}
\begin{align}
    \label{eq:ESS-dynamic}
    & e_{i t}^S = \kappa_i^S e_{i (t-1)}^S - \max \{p_{i t}^S, 0\} \tau / \eta_i^{SD} - \min \{p_{i t}^S, 0\} \tau \eta_i^{SC}, \\
    \label{eq:ESS-p}
    & - \overline{P}_i^{SC} \leq p_{i t}^S \leq \overline{P}_i^{SD}, \\
    \label{eq:ESS-e}
    & \underline{E}_i^S \leq e_{i t}^S \leq \overline{E}_i^S, \\
    \label{eq:ESS-q}
    & (p_{i t}^S)^2 + (q_{i t}^S)^2 \leq (\overline{S}_i^S)^2.
\end{align}
\end{subequations}
In \eqref{eq:ESS-dynamic}, $\kappa_i^S \in (0, 1]$ is a parameter that models the dissipation of ESS $i$. The charging and discharging efficiency coefficients of ESS $i$ are denoted by $\eta_i^{SC} \in (0, 1]$ and $\eta_i^{SD} \in (0, 1]$, respectively. Based on the definition of $p_{i t}^S$, the discharging power is given by $\max \{p_{i t}^S, 0\}$, while the charging power is $- \min \{p_{i t}^S, 0\}$. Therefore, \eqref{eq:ESS-dynamic} captures the SoC dynamics. Constraint \eqref{eq:ESS-p} bounds the power of ESS $i$, where $\overline{P}_i^{SD} > 0$ and $\overline{P}_i^{SC} > 0$ are the maximum discharging and charging power, respectively. Constraint \eqref{eq:ESS-e} defines the lower and upper bounds of the stored energy, while \eqref{eq:ESS-q} specifies the maximum apparent power of ESS $i$.

The following lemma establishes the equivalence between our ESS model and the commonly used model that incorporates a complementarity constraint $p_{i t}^{SD} p_{i t}^{SC} = 0$ with charging power $p_{i t}^{SC}$ and discharging power $p_{i t}^{SD}$ to prevent simultaneous charging and discharging.
\begin{lemma}
\label{lemma:ESS-our-common}
    The set $\{ (e_{i (t - 1)}^S, e_{i t}^S, p_{i t}^S) \in \mathbb{R}^3 \,|\, \eqref{eq:ESS-dynamic},\, \eqref{eq:ESS-p} \}$ equals:
    \begin{align}
    \label{eq:ESS-common}
        \left\{
        \begin{aligned}
            & (e_{i (t - 1)}^S, e_{i t}^S, p_{i t}^S) \in \mathbb{R}^3 \,|\, \exists (p_{i t}^{SD}, p_{i t}^{SC}) \in \mathbb{R}^2,\, \\
            &  \text{s.t.}\; p_{i t}^S = p_{i t}^{SD} - p_{i t}^{SC},\, p_{i t}^{SD} p_{i t}^{SC} = 0, \\
            & 0 \leq p_{i t}^{SD} \leq \overline{P}_i^{SD},\, 0 \leq p_{i t}^{SC} \leq \overline{P}_i^{SC}, \\
            & e_{i t}^S = \kappa_i^S e_{i (t-1)}^S - p_{i t}^{SD} \tau / \eta_i^{SD} + p_{i t}^{SC} \tau \eta_i^{SC}
        \end{aligned}
        \right\}.
    \end{align}
\end{lemma}
The proof of Lemma~\ref{lemma:ESS-our-common} is provided in~\ref{sec:proof}.
We define a useful function to characterize the relationship between the ESS power $p_{i t}^S$ and the change in stored energy in \eqref{eq:ESS-dynamic}:
\begin{definition}
\label{def:es-power}
    For $\eta_i^{SD}, \eta_i^{SC} \in (0, 1]$, define a function $F$ with parameters $\eta_i^{SD}$ and $\eta_i^{SC}$ by:
    \begin{align}
        F_{\eta_i^{SD}, \eta_i^{SC}}(p_{i t}^S) = \max \{p_{i t}^S, 0\} \tau / \eta_i^{SD} + \min \{p_{i t}^S, 0\} \tau \eta_i^{SC}. \nonumber
    \end{align}
\end{definition}

The function $F$ is continuous. Its derivative is $\tau / \eta_i^{SD}$ on $(0, +\infty)$ and $\tau \eta_i^{SC}$ on $(-\infty, 0)$, where $0 < \tau \eta_{i t}^{SC} \leq \tau \leq \tau / \eta_i^{SD} < + \infty$. Consequently, it exhibits the following properties:

\begin{lemma}
\label{lemma:es-power}
    For any $\eta_i^{SD}, \eta_i^{SC} \in (0, 1]$, the function $F_{\eta_i^{SD}, \eta_i^{SC}}: \mathbb{R} \rightarrow \mathbb{R}$ in Definition~\ref{def:es-power} is bijective and strictly increasing. Thus, it has a strictly increasing inverse function $F_{\eta_i^{SD}, \eta_i^{SC}}^{-1}: \mathbb{R} \rightarrow \mathbb{R}$.
\end{lemma}

The ESS model \eqref{eq:ESS} is generally nonconvex due to \eqref{eq:ESS-dynamic}, although all other constraints are convex. A convex approximation was proposed by \cite{shen2020modeling}, which replaces the complementarity constraint $p_{i t}^{SD} p_{i t}^{SC} = 0$ with the linear constraint $p_{i t}^{SD} / \overline{P}_i^{SD} + p_{i t}^{SC} / \overline{P}_i^{SC} \leq 1$. Since the constraints $0 \leq p_{i t}^{SD} \leq \overline{P}_i^{SD},\, 0 \leq p_{i t}^{SC} \leq \overline{P}_i^{SC},\, p_{i t}^{SD} p_{i t}^{SC} = 0$
imply
\begin{align}
    \label{eq:complementary-relax}
    0 \leq p_{i t}^{SD} \leq \overline{P}_i^{SD},\, 0 \leq p_{i t}^{SC} \leq \overline{P}_i^{SC},\, p_{i t}^{SD} / \overline{P}_i^{SD} + p_{i t}^{SC} / \overline{P}_i^{SC} \leq 1,
\end{align}
this approximation relaxes the complementarity constraint to some extent. However, it ensures that ESSs cannot charge and discharge simultaneously in the following sense: Regard $p_{i t}^{SC}$ and $p_{i t}^{SD}$ as the average charging and discharging power of ESS $i$ in period $t$, respectively. Constraint \eqref{eq:complementary-relax} allows ESS $i$ to discharge at power $\overline{P}_i^{SD}$ for a duration of $p_{i t}^{SD} \tau / \overline{P}_i^{SD}$ and then charge at power $\overline{P}_i^{SC}$ for a duration of $p_{i t}^{SC} \tau / \overline{P}_i^{SC}$. The condition $p_{i t}^{SD} / \overline{P}_i^{SD} + p_{i t}^{SC} / \overline{P}_i^{SC} \leq 1$ ensures that the total time required for these operations does not exceed the period length $\tau$, while achieving the desired average power levels $p_{i t}^{SC}$ and $p_{i t}^{SD}$. Consequently, we refer to the following constraints as the \emph{convex ESS constraints}: For any $i \in S_S$ and $t \in S_T$,
\begin{subequations}
\label{eq:ESS-convex}
\begin{align}
    & e_{i t}^S = \kappa_i^S e_{i (t-1)}^S - p_{i t}^{SD} \tau / \eta_i^{SD} + p_{i t}^{SC} \tau \eta_i^{SC}, \\
    & p_{i t}^S = p_{i t}^{SD} - p_{i t}^{SC},\, \eqref{eq:ESS-e},\, \eqref{eq:ESS-q},\, \eqref{eq:complementary-relax}.
\end{align}
\end{subequations}
The original model \eqref{eq:ESS} is referred to as the \emph{general ESS constraints}.

The following lemma demonstrates that for ideal ESSs, the complementarity constraint can be omitted. Consequently, the convex model \eqref{eq:ESS-convex} provides an exact representation for ideal ESSs.
\begin{lemma}
\label{lemma:ESS-ideal}
    If $\eta_i^{SD} = \eta_i^{SC} = 1$, then $\{ (e_{i (t - 1)}^S, e_{i t}^S, p_{i t}^S) \in \mathbb{R}^3 \,|\, \eqref{eq:ESS-dynamic},\, \eqref{eq:ESS-p} \}$ equals the following set:
    \begin{align}
    \label{eq:ESS-relaxed}
        \left\{
        \begin{aligned}
            & (e_{i (t - 1)}^S, e_{i t}^S, p_{i t}^S) \in \mathbb{R}^3 \,|\, \exists (p_{i t}^{SD}, p_{i t}^{SC}) \in \mathbb{R}^2,\, \\
            &  \text{s.t.}\; p_{i t}^S = p_{i t}^{SD} - p_{i t}^{SC},\, 0 \leq p_{i t}^{SD} \leq \overline{P}_i^{SD},\, 0 \leq p_{i t}^{SC} \leq \overline{P}_i^{SC}, \\
            & e_{i t}^S = \kappa_i^S e_{i (t-1)}^S - p_{i t}^{SD} \tau + p_{i t}^{SC} \tau 
        \end{aligned}
        \right\}.
    \end{align}
\end{lemma}
The proof of Lemma~\ref{lemma:ESS-ideal} is provided in~\ref{sec:proof}.

\subsubsection{Power Network}

Let $p_{i t}^N$ and $q_{i t}^N$ represent the net active and reactive power injections at node $i \in S_N$ in period $t \in S_T$, respectively. The aggregate demand $p_t^A$ in period $t$ is the active power injection at node $1$. Therefore, these values are governed by the following nodal power balance equations:
\begin{align}
\label{eq:power-balance}
    p_{i t}^N = p_{i t}^G + p_{i t}^S - p_{i t}^D,\, q_{i t}^N = q_{i t}^G + q_{i t}^S - q_{i t}^D,\, p_t^A = p_{1 t}^N.
\end{align}
If $i \notin S_D$, i.e., there is no load at node $i$, then $p_{i t}^D = q_{i t}^D = 0$ for any $t \in S_T$. Generators and ESSs are treated similarly.

We assume nodal power injections $p_t^N = (p_{i t}^N; i \in S_N)$ and $q_t^N = (q_{i t}^N; i \in S_N)$ are constrained by convex power flow models. For any $t \in S_T$, this is formally expressed as:
\begin{align}
\label{eq:power-flow}
    (p_t^N, q_t^N) \in C_N,
\end{align}
where $C_N$ is a convex set. Common convex power flow models include DC power flow model, linearized DistFlow model \citep{baran1989network}, and SOCP-relaxed branch flow model for radial networks \citep{farivar2013branch}, which are widely used in power system operation.

\subsection{Aggregation Model}

By integrating the component models in Section~\ref{sec:model-component}, we formulate the operational constraints of the distribution system: For any $t \in S_T$,
\begin{subequations}
\label{eq:cons}
\begin{align}
    \label{eq:cons-1}
    & (p_t^A, p_t^S, y_t) \in C_t, \\
    \label{eq:cons-2}
    & e_{i t}^S = \kappa_i^S e_{i (t-1)}^S - \max \{p_{i t}^S, 0\} \tau / \eta_i^{SD} - \min \{p_{i t}^S, 0\} \tau \eta_i^{SC},\, \forall i \in S_S, \\
    \label{eq:cons-3}
    & \underline{E}_i^S \leq e_{i t}^S \leq \overline{E}_i^S,\, \forall i \in S_S.
\end{align}
\end{subequations}
For conciseness, time-decoupled constraints \eqref{eq:load}, \eqref{eq:DER}, \eqref{eq:ESS-p}, \eqref{eq:ESS-q}, \eqref{eq:power-balance}, and \eqref{eq:power-flow} are collected into \eqref{eq:cons-1}, where $p_t^S = (p_{i t}^S; i \in S_S)$ is the ESS active power vector; $y_t = (p_{i t}^D, q_{i t}^D, p_{i t}^G, q_{i t}^G, q_{i t}^S, p_{i t}^N, q_{i t}^N; i \in S_N)$ stacks other variables; $C_t$ is a convex and closed set encoding the constraints in period $t$. The general ESS model \eqref{eq:ESS} is used in \eqref{eq:cons}. The operational constraints using the convex ESS model \eqref{eq:ESS-convex} can be similarly written as: For any $t \in S_T$,
\begin{subequations}
\label{eq:cons-convex}
\begin{align}
    \label{eq:cons-convex-1}
    & (p_t^A, p_t^{SD}, p_t^{SC}, y_t) \in C_t', \\
    \label{eq:cons-convex-2}
    & e_{i t}^S = \kappa_i^S e_{i (t-1)}^S - p_{i t}^{SD} \tau / \eta_i^{SD} + p_{i t}^{SC} \tau \eta_i^{SC},\, \forall i \in S_S, \\
    \label{eq:cons-convex-3}
    & \underline{E}_i^S \leq e_{i t}^S \leq \overline{E}_i^S,\, \forall i \in S_S,
\end{align}
\end{subequations}
where the convex and closed set $C_t'$ collects constraints \eqref{eq:load}, \eqref{eq:DER}, \eqref{eq:ESS-q}, \eqref{eq:complementary-relax}, \eqref{eq:power-balance}, and \eqref{eq:power-flow}.
In the following, we first introduce the two-stage aggregation model, followed by the multistage model.

\subsubsection{Two-Stage Aggregation Model}

In time-decoupled power flexibility aggregation, the goal is to find the ``largest" set $S_A = \times_{t \in S_T} [p_t^{A \vee}, p_t^{A \wedge}]$ such that any aggregate power trajectory $p^A = (p_t^A; t \in S_T) \in S_A$ can be realized through distribution system operation management. The objective function is denoted by $\phi(p^{A \vee}, p^{A \wedge})$, where $\phi$ is a flexibility index function. A common choice is $\phi(p^{A \vee}, p^{A \wedge}) = \sum_{t \in S_T} \omega_t (p_t^{A \wedge} - p_t^{A \vee})$, where $\omega = (\omega_t; t \in S_T)$ is a weight vector indicating the importance of flexibility in each period, with $\omega_t \geq 0$ for all $t \in S_T$. When $\omega_t = 1$ for all $t \in S_T$, all time periods are equally weighted \citep{chen2021leveraging}. Another choice is $\phi(p^{A \vee}, p^{A \wedge}) = \sum_{t \in S_T} \log (p_t^{A \wedge} - p_t^{A \vee})$, while maximizing it is equivalent to maximizing the volume $\Pi_{t \in S_T} (p_t^{A \wedge} - p_t^{A \vee})$ of $S_A$. Based on the objective function $\phi(p^{A \vee}, p^{A \wedge})$, the \emph{Two-Stage} model for aggregation is as follows:
\begin{subequations}
\label{eq:two-stage}
\begin{align}
    \label{eq:two-stage-obj}
    \textbf{Two-Stage}:\; \max_{p^{A \vee}, p^{A \wedge}}\; & \phi(p^{A \vee}, p^{A \wedge}) \\
    \label{eq:two-stage-con}
    \text{s.t.}\; & \forall p^A \in \times_{t \in S_T} [p_t^{A \vee}, p_t^{A \wedge}],\; \exists e^S, p^S, y,\, \text{s.t.}\; \eqref{eq:cons},\, \forall t \in S_T.
\end{align}
\end{subequations}
In \eqref{eq:two-stage}, the first-stage decision variables are $p^{A \wedge}$ and $p^{A \vee}$, which define the time-decoupled aggregate power flexibility region. The aggregate power $p^A$ (setpoints by TSO) is the uncertainty in this two-stage RO model, which is decision-dependent since its range is determined by the first-stage decision variables. Constraint \eqref{eq:two-stage-con} ensures that for any $p^A$ satisfying $p_t^{A \vee} \leq p_t^A \leq p_t^{A \wedge}$ for all $t \in S_T$, there exists an operation strategy characterized by the second-stage decision variables $e^S$, $p^S$, and $y$ such that the operational constraints \eqref{eq:cons} are satisfied. However, since the operation strategy is determined only after the entire aggregate power trajectory is known, model \eqref{eq:two-stage} is two-stage and cannot account for the TSO's sequential decision-making process.

\subsubsection{Multistage Aggregation Model}

The proposed \emph{Multistage} model for aggregation is as follows:
\begin{subequations}
\label{eq:multistage}
\begin{align}
    \textbf{Multistage}:\; \max_{p^{A \vee}, p^{A \wedge}}\; & \phi(p^{A \vee}, p^{A \wedge}) \\
    \text{s.t.}\; & \forall p_1^A \in [p_1^{A \vee}, p_1^{A \wedge}],\, \exists e_1^S, p_1^S, y_1,\, \text{s.t.} \nonumber \\
    & \forall p_2^A \in [p_2^{A \vee}, p_2^{A \wedge}],\, \exists e_2^S, p_2^S, y_2,\, \text{s.t.} \nonumber \\
    & \dots \nonumber \\
    & \forall p_T^A \in [p_T^{A \vee}, p_T^{A \wedge}],\, \exists e_T^S, p_T^S, y_T,\, \text{s.t.}\; \eqref{eq:cons},\, \forall t \in S_T.
\end{align}
\end{subequations}
In the Multistage model \eqref{eq:multistage}, the decisions of the DSO and the TSO are made sequentially. In the first stage, the DSO determines the bounds $p^{A \wedge}$ and $p^{A \vee}$. Subsequently, the TSO reveals the aggregate power (setpoint) $p_1^A$ and the DSO determines the operation strategies $e_1^S$, $p_1^S$, and $y_1$ without knowledge of future setpoints $p_2^A, p_3^A, \dots, p_T^A$. This process continues and for each period $t \in S_T$, the DSO determines the operation strategies $e_t^S$, $p_t^S$, and $y_t$ based solely on the revealed uncertainties $p_1^A, p_2^A, \dots, p_t^A$. In this way, model \eqref{eq:multistage} accounts for the TSO's sequential actions and the corresponding nonanticipativity constraint. Since the range of the uncertainty $p^A$ depends on the first-stage decision variables $p^{A \wedge}$ and $p^{A \vee}$, model \eqref{eq:multistage} is a multistage RO problem with DDU.

The Multistage model \eqref{eq:multistage} is similar to the Two-Stage model \eqref{eq:two-stage} except that the operation strategies in each period must be independent of future realizations of the aggregate power. This additional requirement makes the Multistage model more restrictive than the Two-Stage model. Consequently, the Two-Stage model can be viewed as a relaxation of the Multistage model. These are summarized the following proposition:
\begin{proposition}
\label{prop:two-stage}
    The feasible region of the Multistage model \eqref{eq:multistage} is a subset of the feasible region of the Two-Stage model \eqref{eq:two-stage}. Furthermore, their optimal values satisfy $\eqref{eq:multistage} \leq \eqref{eq:two-stage}$.
\end{proposition}

Model \eqref{eq:multistage} is referred to as the \emph{Multistage model under general ESS constraints}, because the general ESS model \eqref{eq:ESS} is incorporated in \eqref{eq:cons}. We also investigate the \emph{Multistage model under convex ESS constraints}, which is analogous to \eqref{eq:multistage} but replaces constraint \eqref{eq:cons} with the convex formulation \eqref{eq:cons-convex}:
\begin{subequations}
\label{eq:multistage-convex}
\begin{align}
    & \textbf{Multistage (Convex ESS)}: \nonumber \\
    \max_{p^{A \vee}, p^{A \wedge}}\; & \phi(p^{A \vee}, p^{A \wedge}) \\
    \text{s.t.}\; & \forall p_1^A \in [p_1^{A \vee}, p_1^{A \wedge}],\, \exists e_1^S, p_1^{SD}, p_1^{SC}, y_1,\, \text{s.t.} \nonumber \\
    & \forall p_2^A \in [p_2^{A \vee}, p_2^{A \wedge}],\, \exists e_2^S, p_2^{SD}, p_2^{SC}, y_2,\, \text{s.t.} \nonumber \\
    & \dots \nonumber \\
    & \forall p_T^A \in [p_T^{A \vee}, p_T^{A \wedge}],\, \exists e_T^S, p_T^{SD}, p_T^{SC}, y_T,\, \text{s.t.}\; \eqref{eq:cons-convex},\, \forall t \in S_T.
\end{align}
\end{subequations}

\section{Solution Methods}
\label{sec:solution}

In this section, we present four solution methods for multistage time-decoupled power flexibility aggregation: An exact solution method based on enumeration is proposed for the Multistage model under convex ESS constraints. For general ESS constraints, an inner approximation solution method based on SoC ranges is proposed, and the envelope-based solution method \citep{chen2019aggregate} is analyzed. For these exact and inner approximation methods, greedy disaggregation algorithms are established considering operation costs. Furthermore, we derive outer approximations through the Two-Stage model and its relaxation. Finally, we establish theoretical comparisons of these methods in terms of applicability and accuracy.

\subsection{Exact Solution Based on Enumeration}

This subsection develops an exact solution for the Multistage model under convex ESS constraints. We establish the aggregation and disaggregation methods in sequence.

\subsubsection{Aggregation}

ESSs introduce temporal coupling in the distribution system through their SoC dynamics. To facilitate analysis, we define the feasible set of SoC ranges as follows:
\begin{definition}
\label{def:SoC-range}
    Suppose $(p^{A \vee}, p^{A \wedge})$ is a solution to the Multistage model \eqref{eq:multistage}. For $t_0 \in S_T \cup \{0\}$, define set
    $\mathcal{E}_{t_0}^S(p^{A \vee}, p^{A \wedge})$ as:
\begin{align}
\label{eq:set-soc}
    \mathcal{E}_{t_0}^S(p^{A \vee}, p^{A \wedge}) = \left\{ e_{t_0}^S \,\middle|\,  
    \begin{aligned}
        & \underline{E}_i^S \leq e_{i t_0}^S \leq \overline{E}_i^S,\, \forall i \in S_S, \\
        & \forall p_{t_0 + 1}^A \in [p_{t_0 + 1}^{A \vee}, p_{t_0 + 1}^{A \wedge}],\, \exists e_{t_0 + 1}^S, p_{t_0 + 1}^S, y_{t_0 + 1},\, \text{s.t.} \\
        & \forall p_{t_0 + 2}^A \in [p_{t_0 + 2}^{A \vee}, p_{t_0 + 2}^{A \wedge}],\, \exists e_{t_0 + 2}^S, p_{t_0 + 2}^S, y_{t_0 + 2},\, \text{s.t.} \\
        & \dots \\
        & \forall p_T^A \in [p_T^{A \vee}, p_T^{A \wedge}],\, \exists e_T^S, p_T^S, y_T,\, \text{s.t.}\; \eqref{eq:cons},\, \forall t > t_0
    \end{aligned}
    \right\}.
\end{align}
\end{definition}

Given $(p^{A \vee}, p^{A \wedge})$ and $t_0 \in S_T \cup \{0\}$, $\mathcal{E}_{t_0}^S(p^{A \vee}, p^{A \wedge})$ is the set of feasible SoC values of ESSs at the end of period $t_0$, considering the operational constraints after period $t_0$ and the nonanticipativity constraints. According to Definition~\ref{def:SoC-range}, $\mathcal{E}_{t_0}^S(p^{A \vee}, p^{A \wedge})$ has the following properties:

\begin{lemma}
    \label{lemma:set-soc}
    Suppose $(p^{A \vee}, p^{A \wedge})$ is a solution to the Multistage model \eqref{eq:multistage}. Then $(p^{A \vee}, p^{A \wedge})$ is feasible in \eqref{eq:multistage} if and only if the initial SoC vector $e_0^S \in \mathcal{E}_0^S(p^{A \vee}, p^{A \wedge})$. In this case, for any $t_0 \in S_T \cup \{0\}$, $\mathcal{E}_{t_0}^S(p^{A \vee}, p^{A \wedge})$ is a nonempty and closed subset of $\{ e_{t_0}^S \,|\, \underline{E}_i^S \leq e_{i t_0}^S \leq \overline{E}_i^S, \forall i \in S_S \}$.
\end{lemma}

The following proposition addresses the case of convex ESS constraints, showing that the interval $[p_t^{A \vee}, p_t^{A \wedge}]$ can be relaxed to $\{p_t^{A \vee}, p_t^{A \wedge}\}$ for $t > t_0$ without altering the set $\mathcal{E}_{t_0}^S(p^{A \vee}, p^{A \wedge})$.

\begin{proposition}
\label{prop:set-soc-reformulate}
    Consider the Multistage model \eqref{eq:multistage-convex} under convex ESS constraints. Suppose $(p^{A \vee}, p^{A \wedge})$ is feasible in \eqref{eq:multistage-convex}. Then for any $t_0 \in S_T \cup \{0\}$, $\mathcal{E}_{t_0}^S(p^{A \vee}, p^{A \wedge})$ equals:
    \begin{align}
    \label{eq:set-soc-reformulate}
        \left\{ e_{t_0}^S \,\middle|\,  
        \begin{aligned}
            & \underline{E}_i^S \leq e_{i t_0}^S \leq \overline{E}_i^S,\, \forall i \in S_S, \\
            & \forall p_{t_0 + 1}^A \in \{p_{t_0 + 1}^{A \vee}, p_{t_0 + 1}^{A \wedge}\},\, \exists e_{t_0 + 1}^S, p_{t_0 + 1}^{SD}, p_{t_0 + 1}^{SC}, y_{t_0 + 1},\, \text{s.t.} \\
            & \forall p_{t_0 + 2}^A \in \{p_{t_0 + 2}^{A \vee}, p_{t_0 + 2}^{A \wedge}\},\, \exists e_{t_0 + 2}^S, p_{t_0 + 2}^{SD}, p_{t_0 + 2}^{SC}, y_{t_0 + 2},\, \text{s.t.} \\
            & \dots \\
            & \forall p_T^A \in \{p_T^{A \vee}, p_T^{A \wedge}\},\, \exists e_T^S, p_T^{SD}, p_T^{SC}, y_T,\, \text{s.t.}\; \eqref{eq:cons-convex},\, \forall t > t_0
        \end{aligned}
        \right\}.
    \end{align}
\end{proposition}

The proof of Proposition~\ref{prop:set-soc-reformulate} is provided in \ref{sec:proof}. Because the constraint of the Multistage model is equivalent to $e_0^S \in \mathcal{E}_0^S(p^{A \vee}, p^{A \wedge})$, we have the following corollary:
\begin{corollary}
\label{cor:enumeration}
    The Multistage model \eqref{eq:multistage-convex} under convex ESS constraints is equivalent to:
    \begin{subequations}
    \label{eq:multistage-enumeration}
    \begin{align}
        \textbf{Enumeration}:\; \max_{p^{A \vee}, p^{A \wedge}}~ & \phi(p^{A \vee}, p^{A \wedge}) \\
        \label{eq:pa-order}
        \text{s.t.}~ & p_t^{A \vee} \leq p_t^{A \wedge},\, \forall t \in S_T, \\
        & \forall p^A \in \times_{t \in S_T} \{p_t^{A \vee}, p_t^{A \wedge}\}, \nonumber \\
        & \exists e_t^S(p_{\leq t}^A), p_t^{SD} (p_{\leq t}^A), p_t^{SC} (p_{\leq t}^A), y_t (p_{\leq t}^A); t \in S_T,\, \nonumber \\
        & \text{s.t.}\; \eqref{eq:cons-convex},\, \forall t \in S_T.
    \end{align}
    \end{subequations}
\end{corollary}

In \eqref{eq:multistage-enumeration}, the uncertainty set $\times_{t \in S_T} \{p_t^{A \vee}, p_t^{A \wedge}\}$ contains at most $2^T$ distinct elements, which can be effectively enumerated if $T$ is not large. This enables us to treat \eqref{eq:multistage-enumeration} as a single-stage optimization problem. The notation $e_t^S(p_{\leq t}^A)$ indicates that the value of $e_t^S$ depends on $p_{\leq t}^A = (p_1^A, p_2^A, \dots, p_t^A)$. Since there are $2^t$ possible values of $p_{\leq t}^A$, $e_t^S(p_{\leq t}^A)$ can be interpreted as $2^t$ copies of the variable $e_t^S$. If $\tilde{p}_{\leq t_0}^A = \hat{p}_{\leq t_0}^A$, then $e_t^S(\tilde{p}_{\leq t}^A) = e_t^S(\hat{p}_{\leq t}^A)$ for any $t \leq t_0$. The notations $p_t^{SD} (p_{\leq t}^A)$, $p_t^{SC} (p_{\leq t}^A)$, and $y_t (p_{\leq t}^A)$ are similar. Due to the convexity of \eqref{eq:cons-convex}, \eqref{eq:multistage-enumeration} becomes a single-stage convex program after the enumeration if $-\phi(p^{A \vee}, p^{A \wedge})$ is a convex function.

Model \eqref{eq:multistage-enumeration} provides an exact solution for the Multistage model \eqref{eq:multistage-convex} under convex ESS constraints by enumerating the elements of the uncertainty set and solving a single-stage program. We refer to it as the \emph{Enumeration} model. However, the number of variables and constraints grows exponentially with $T$, making this method impractical for cases with many periods.

\subsubsection{Disaggregation}

Consider the Multistage model \eqref{eq:multistage-convex} under convex ESS constraints. While the proof of Proposition~\ref{prop:set-soc-reformulate} provides a method for disaggregation, it does not incorporate the DSO's objective of minimizing operation costs. Therefore, we propose a greedy disaggregation algorithm for the Enumeration method, which minimizes the operation cost in each period. Suppose the aggregation result sent by the DSO to the TSO is represented by $(p^{A \vee}, p^{A \wedge})$, a feasible solution to model \eqref{eq:multistage-convex}. The TSO sequentially determines the aggregate power trajectory $p^A$, ensuring $p_t^A \in [p_t^{A \vee}, p_t^{A \wedge}]$ for all $t \in S_T$. For a period $t_0 \in S_T$, all variables before $t_0$ have been determined. The aggregate power $p_{t_0}^A$ is specified by the TSO, while the future aggregate power $p_{> t_0}^A$ remains unknown. Our task is to determine $p_{t_0}^{SD}$, $p_{t_0}^{SC}$, and $y_{t_0}$ to realize $p_{t_0}^A$, ensure future feasibility, and minimize the operation cost $c(p_{t_0}^{SD}, p_{t_0}^{SC}, y_{t_0})$ for period $t_0$. Similar to Proposition~\ref{prop:set-soc-reformulate}, the uncertainty set $\times_{t_0 < t \in S_T} [p_t^{A \vee}, p_t^{A \wedge}]$ of future aggregate power can be replaced by $\times_{t_0 < t \in S_T} \{p_t^{A \vee}, p_t^{A \wedge}\}$. To this end, we develop the following program for real-time disaggregation: For $t_0 \in S_T$,
\begin{subequations}
\label{eq:multistage-enumeration-greedy}
\begin{align}
    \min_{p_{t_0}^{SD}, p_{t_0}^{SC}, y_{t_0}}\; & c(p_{t_0}^{SD}, p_{t_0}^{SC}, y_{t_0}) \\
    \text{s.t.}\; & \eqref{eq:cons-convex},\, \text{for}\; t = t_0, \\
    & \forall p_{> t_0}^A \in \times_{t_0 < t \in S_T} \{p_t^{A \vee}, p_t^{A \wedge}\}, \nonumber \\
    & \exists e_t^S(p_{\leq t}^A), p_t^{SD} (p_{\leq t}^A), p_t^{SC} (p_{\leq t}^A), y_t (p_{\leq t}^A); t_0 < t \in S_T,\, \nonumber \\
    & \text{s.t.}\; \eqref{eq:cons-convex},\, \forall t_0 < t \in S_T.
\end{align}
\end{subequations}
Similar to \eqref{eq:multistage-enumeration}, model \eqref{eq:multistage-enumeration-greedy} reduces to a single-stage program after enumerating the $2^{T - t_0}$ elements in the uncertainty set $\times_{t_0 < t \in S_T} \{p_t^{A \vee}, p_t^{A \wedge}\}$. Moreover, it is convex if the operation cost function $c(p_{t_0}^{SD}, p_{t_0}^{SC}, y_{t_0})$ is convex. We summarize the aggregation and disaggregation processes of the Enumeration method in Algorithm~\ref{alg:enumeration}.

\begin{algorithm}
    \caption{Aggregation and disaggregation of the Enumeration method}
    \begin{algorithmic}[1]
    \label{alg:enumeration}
        \STATE \textbf{DSO:}
        \STATE \quad Aggregation: Solve \eqref{eq:multistage-enumeration} to obtain $p^{A \vee}$ and $p^{A \wedge}$. 
        \STATE \quad Send $p^{A \vee}$ and $p^{A \wedge}$ to the TSO.
        \FOR{$t_0 = 1$ \TO $T$}
        \STATE \textbf{TSO:}
        \STATE \quad Determine $p_{t_0}^A \in [p_{t_0}^{A \vee}, p_{t_0}^{A \wedge}]$. Send $p_{t_0}^A$ to the DSO.
        \STATE \textbf{DSO:}
        \STATE \quad Disaggregation: Solve \eqref{eq:multistage-enumeration-greedy} to obtain $e_{t_0}^S$, $p_{t_0}^{SD}$, $p_{t_0}^{SC}$, and $y_{t_0}$.
        \ENDFOR
    \end{algorithmic}
\end{algorithm}

\subsection{Inner Approximation Based on SoC Ranges}

This subsection develops an inner approximation method for the Multistage model \eqref{eq:multistage} under general ESS constraints. First, we derive an equivalent formulation based on SoC ranges. Next, we propose an inner approximation using rectangular SoC ranges. Finally, we introduce the corresponding disaggregation method.

\subsubsection{Equivalent Formulation Based on SoC Ranges}

SoC ranges play a critical role in the Multistage model \eqref{eq:multistage} because the SoC dynamics constraint in \eqref{eq:ESS-dynamic} is the only time-coupled constraint. Recall the set $\mathcal{E}_t^S(p^{A \vee}, p^{A \wedge})$ defined in Definition~\ref{def:SoC-range} for $t \in S_T$, which represents the feasible SoC range at the end of period $t$. We can express $\mathcal{E}_{t - 1}^S(p^{A \vee}, p^{A \wedge})$ in terms of $\mathcal{E}_t^S(p^{A \vee}, p^{A \wedge})$ as the following lemma shows:

\begin{lemma}
\label{lemma:soc-range}
    Suppose $(p^{A \vee}, p^{A \wedge})$ is a solution to the Multistage model \eqref{eq:multistage}. For $t \in S_T$, let the sets
    $\mathcal{E}_t^S(p^{A \vee}, p^{A \wedge})$ and $\mathcal{E}_{t - 1}^S(p^{A \vee}, p^{A \wedge})$ be defined as in Definition~\ref{def:SoC-range}. Then $\mathcal{E}_{t - 1}^S(p^{A \vee}, p^{A \wedge})$ equals:
    \begin{align}
        \left\{ e_{t - 1}^S \,\middle|\,  
    \begin{aligned}
        & \underline{E}_i^S \leq e_{i (t - 1)}^S \leq \overline{E}_i^S,\, \forall i \in S_S, \\
        & \forall p_t^A \in [p_t^{A \vee}, p_t^{A \wedge}],\, \exists e_t^S, p_t^S, y_t,\, \text{s.t.}\; e_t^S \in \mathcal{E}_t^S(p^{A \vee}, p^{A \wedge}),\, \eqref{eq:cons}
    \end{aligned}
    \right\}. \nonumber
    \end{align}
\end{lemma}

Using Lemma~\ref{lemma:soc-range}, we can reformulate the constraint of the Multistage model \eqref{eq:multistage} in terms of SoC ranges for adjacent periods, as shown in the following proposition:

\begin{proposition}
\label{prop:set-soc-equivalent}
    Suppose $(p^{A \vee}, p^{A \wedge})$ is a solution to the Multistage model \eqref{eq:multistage}. Then $(p^{A \vee}, p^{A \wedge})$ is feasible in \eqref{eq:multistage} if and only if there are nonempty and closed sets $\mathcal{E}_t^S \subseteq \{ e_t^S \,|\, \underline{E}_i^S \leq e_{i t}^S \leq \overline{E}_i^S, \forall i \in S_S \}$ for all $t \in S_T \cup \{0\}$ such that $e_0^S \in \mathcal{E}_0^S$ and for any $t \in S_T$,
    \begin{align}
    \label{eq:set-soc-equivalent}
        \forall e_{t - 1}^S \in \mathcal{E}_{t - 1}^S,\, \forall p_t^A \in [p_t^{A \vee}, p_t^{A \wedge}],\, \exists e_t^S \in \mathcal{E}_t^S, p_t^S, y_t, \, \text{s.t.}\; \eqref{eq:cons-1}, \eqref{eq:cons-2}.
    \end{align}
\end{proposition}

The proof of Proposition~\ref{prop:set-soc-equivalent} is provided in \ref{sec:proof}. Using Proposition~\ref{prop:set-soc-equivalent}, we derive the following equivalent formulation of the Multistage model \eqref{eq:multistage}:

\begin{corollary}
\label{cor:multi-stage-equivalent}
    The Multistage model \eqref{eq:multistage} has the following equivalent formulation:
    \begin{subequations}
    \label{eq:multistage-equivalent}
    \begin{align}
        \max_{p_t^{A \vee}, p_t^{A \wedge}, \mathcal{E}_t^S, \forall t}~ & \phi (p^{A \vee}, p^{A \wedge}) \\
        \text{s.t.}~ & e_0^S \in \mathcal{E}_0^S,\, \eqref{eq:set-soc-equivalent},\, \forall t \in S_T, \\
        & \mathcal{E}_t^S \subseteq \{ e_{i t}^S ~|~ \underline{E}_i^S \leq e_{i t}^S \leq \overline{E}_i^S, \forall i \},\, \forall t \in S_T \cup \{0\}.
    \end{align}
    \end{subequations}
\end{corollary}

From the proof of Proposition~\ref{prop:set-soc-equivalent}, it follows that if $(p^{A \vee}, p^{A \wedge})$ is optimal in the Multistage model \eqref{eq:multistage}, then $(p^{A \vee}, p^{A \wedge})$ and $\mathcal{E}_t(p^{A \vee}, p^{A \wedge})$ for $t \in S_T \cup \{0\}$ form an optimal solution to \eqref{eq:multistage-equivalent}. Although the two formulations are equivalent, \eqref{eq:multistage-equivalent} cannot be solved directly because its variables include sets $\mathcal{E}_t^S$ for $t \in S_T \cup \{0\}$. To address this, we will restrict these sets to rectangular forms and derive an inner approximate solution. Before proceeding, we present a result about the convexity of $\mathcal{E}_t^S(p^{A \vee}, p^{A \wedge})$ for feasible $(p^{A \vee}, p^{A \wedge})$ in the Multistage model \eqref{eq:multistage-convex} under convex ESS constraints, whose proof is provided in \ref{sec:proof}:

\begin{proposition}
\label{prop:set-soc-convex}
    Suppose $(p^{A \vee}, p^{A \wedge})$ is feasible in the Multistage model \eqref{eq:multistage-convex} under convex ESS constraints. Then for any $t \in S_T \cup \{0\}$, the set $\mathcal{E}_t^S(p^{A \vee}, p^{A \wedge})$ is convex.
\end{proposition}


\subsubsection{Inner Approximation Based on Rectangular SoC Ranges}

We restrict the set $\mathcal{E}_t^S$ in Proposition~\ref{prop:set-soc-equivalent} and Corollary~\ref{cor:multi-stage-equivalent} to rectangular SoC ranges as follows:
\begin{align}
    \mathcal{E}_t^S = \times_{i \in S_S} \left[ e_{i t}^{S \downarrow}, e_{i t}^{S \uparrow} \right] = \left\{ e_t^S ~\middle|~ e_{i t}^{S \downarrow} \leq e_{i t}^S \leq e_{i t}^{S \uparrow}, \forall i \in S_S \right\}. \nonumber
\end{align}
This yields an inner approximation for the Multistage model \eqref{eq:multistage}:
\begin{subequations}
\label{eq:multistage-box}
\begin{align}
    \max_{p^{A \vee}, p^{A \wedge}, e^{S \downarrow}, e^{S \uparrow}}\; & \phi(p^{A \vee}, p^{A \wedge}) \\
    \label{eq:multistage-box-2}
    \text{s.t.}\; & e_{i 0}^{S \downarrow} \leq e_{i 0}^S \leq e_{i 0}^{S \uparrow},\, \forall i \in S_S, \\
    \label{eq:multistage-box-3}
    & \underline{E}_i^S \leq e_{i t}^{S \downarrow} \leq e_{i t}^{S \uparrow} \leq \overline{E}_i^S,\, \forall i \in S_S,\, \forall t \in S_T \cup \{0\}, \\
    & \forall t \in S_T,\, \forall e_{t - 1}^S \in \times_{i \in S_S} \left[ e_{i (t - 1)}^{S \downarrow}, e_{i (t - 1)}^{S \uparrow} \right],\, \forall p_t^A \in [p_t^{A \vee}, p_t^{A \wedge}], \nonumber \\
    \label{eq:multistage-box-4}
    & \exists e_t^S, p_t^S, y_t,\, \text{s.t.}\; e_{i t}^{S \downarrow} \leq e_{i t}^S \leq e_{i t}^{S \uparrow},\, \forall i \in S_S,\, \eqref{eq:cons-1}, \eqref{eq:cons-2}.
\end{align}
\end{subequations}

The following proposition demonstrates that it suffices to check $\forall p_t^A \in \{ p_t^{A \vee}, p_t^{A \wedge}\}$ rather than $\forall p_t^A \in [p_t^{A \vee}, p_t^{A \wedge}]$ in \eqref{eq:multistage-box-4}.

\begin{proposition}
\label{prop:multistage-box-uncertainty}
    Suppose $p_t^{A \vee} \leq p_t^{A \wedge}$ for any $t \in S_T$. Then \eqref{eq:multistage-box-4} is equivalent to: 
    \begin{subequations}
    \label{eq:set-soc-equivalent-extreme}
    \begin{align}
        & \forall t \in S_T,\, \forall e_{t - 1}^S \in \times_{i \in S_S} \left[ e_{i (t - 1)}^{S \downarrow}, e_{i (t - 1)}^{S \uparrow} \right], \nonumber \\
        & \exists e_t^{S \vee}, p_t^{S \vee}, y_t^\vee, e_t^{S \wedge}, p_t^{S \wedge}, y_t^\wedge,\, \text{s.t.} \nonumber \\
        \label{eq:set-soc-equivalent-extreme-1}
        &  (p_t^{A \vee}, p_t^{S \vee}, y_t^\vee), (p_t^{A \wedge}, p_t^{S \wedge}, y_t^\wedge) \in C_t, \\
        \label{eq:set-soc-equivalent-extreme-2}
        & e_{i t}^{S \downarrow} \leq e_{i t}^{S \vee} \leq e_{i t}^{S \uparrow},\, e_{i t}^{S \downarrow} \leq e_{i t}^{S \wedge} \leq e_{i t}^{S \uparrow},\, \forall i \in S_S, \\
        \label{eq:set-soc-equivalent-extreme-3}
        & e_{i t}^{S \vee} = \kappa_i^S e_{i (t - 1)}^S - F_{\eta_i^{SD}, \eta_i^{SC}} (p_{i t}^{S \vee}),\, \forall i \in S_S, \\
        \label{eq:set-soc-equivalent-extreme-4}
        & e_{i t}^{S \wedge} = \kappa_i^S e_{i (t - 1)}^S - F_{\eta_i^{SD}, \eta_i^{SC}} (p_{i t}^{S \wedge}),\, \forall i \in S_S.
    \end{align}
    \end{subequations}
\end{proposition}

The proof of Proposition~\ref{prop:multistage-box-uncertainty} is provided in \ref{sec:proof}. Based on this result, we propose an inner approximation solution method for the Multistage model \eqref{eq:multistage}, achieved by constraining the SoC ranges to rectangular sets:

\begin{corollary}
\label{cor:rectangular}
    If $(p^{A \vee}, p^{A \wedge}, e^{S \downarrow}, e^{S \uparrow})$ is feasible in the following problem \eqref{eq:multistage-box-two}, then $(p^{A \vee}, p^{A \wedge})$ is feasible in the Multistage model \eqref{eq:multistage}. Furthermore, the optimal values satisfy $\eqref{eq:multistage-box-two} \leq \eqref{eq:multistage}$. 
    \begin{subequations}
    \label{eq:multistage-box-two}
    \begin{align}
        \textbf{Rectangular}:\; \max_{p^{A \vee}, p^{A \wedge}, e^{S \downarrow}, e^{S \uparrow}}\; & \phi(p^{A \vee}, p^{A \wedge}) \\
        \text{s.t.}\; & \eqref{eq:pa-order}, \eqref{eq:multistage-box-2}, \eqref{eq:multistage-box-3}, \eqref{eq:set-soc-equivalent-extreme}. 
    \end{align}
    \end{subequations}
\end{corollary}

The \emph{Rectangular} model \eqref{eq:multistage-box-two} is a two-stage RO problem, which has better tractability than the Multistage model \eqref{eq:multistage}. When the convex ESS constraints are employed and operational constraints are linear, \eqref{eq:multistage-box-two} admits an exact solution through an improved C\&CG algorithm \citep{zeng2022two}. In fact, for the Multistage model \eqref{eq:multistage-convex} under convex ESS constraints, uncertainty sets can be reduced to their vertex sets due to the convexity of the second-stage constraints, and thus the rectangular inner approximation is as follows:
\begin{subequations}
\label{eq:multistage-box-two-enumeration}
\begin{align}
    \max_{p^{A \vee}, p^{A \wedge}, e^{S \downarrow}, e^{S \uparrow}}\; & \phi(p^{A \vee}, p^{A \wedge}) \\
    \text{s.t.}\; & \eqref{eq:pa-order}, \eqref{eq:multistage-box-2}, \eqref{eq:multistage-box-3}, \\
    & \forall t \in S_T,\, \forall e_{t - 1}^S \in \times_{i \in S_S} \left\{ e_{i (t - 1)}^{S \downarrow}, e_{i (t - 1)}^{S \uparrow} \right\}, \nonumber \\
    & \exists e_t^{S \vee}, p_t^{SD \vee}, p_t^{SC \vee}, y_t^\vee, e_t^{S \wedge}, p_t^{SD \wedge}, p_t^{SC \wedge}, y_t^\wedge,\, \text{s.t.} \nonumber \\
    & (p_t^{A \vee}, p_t^{SD \vee}, p_t^{SC \vee}, y_t^\vee), (p_t^{A \wedge}, p_t^{SD \wedge}, p_t^{SC \wedge}, y_t^\wedge) \in C_t, \nonumber \\
    & e_{i t}^{S \downarrow} \leq e_{i t}^{S \vee} \leq e_{i t}^{S \uparrow},\, e_{i t}^{S \downarrow} \leq e_{i t}^{S \wedge} \leq e_{i t}^{S \uparrow},\, \forall i \in S_S, \nonumber \\
    & e_{i t}^{S \vee} = \kappa_i^S e_{i (t - 1)}^S - p_t^{SD \vee} \tau / \eta_i^{SD} + p_t^{SC \vee} \tau \eta_i^{SC},\, \forall i \in S_S, \nonumber \\
    \label{eq:multistage-box-two-enumeration-4}
    & e_{i t}^{S \wedge} = \kappa_i^S e_{i (t - 1)}^S - p_t^{SD \wedge} \tau / \eta_i^{SD} + p_t^{SC \wedge} \tau \eta_i^{SC},\, \forall i \in S_S.
\end{align}
\end{subequations}
When the number of ESSs is small, the vertices can be enumerated, allowing \eqref{eq:multistage-box-two-enumeration} to be solved as a single-stage program.

We continue to discuss the case with a single ESS in the distribution system, where we find that the Rectangular model \eqref{eq:multistage-box-two} provides an exact solution to the Multistage model \eqref{eq:multistage}, considering general ESS constraints. To show this, we first prove the convexity of $\mathcal{E}_t^S(p^{A \vee}, p^{A \wedge})$ when there is only one ESS.

\begin{proposition}
\label{prop:set-soc-convex-1es}
    Suppose there is only one ESS and $(p^{A \vee}, p^{A \wedge})$ is feasible in the Multistage model \eqref{eq:multistage}. Then for any $t \in S_T \cup \{0\}$, the set $\mathcal{E}_t^S(p^{A \vee}, p^{A \wedge})$ in Definition~\ref{def:SoC-range} is convex.
\end{proposition}

The proof of Proposition~\ref{prop:set-soc-convex-1es} is provided in \ref{sec:proof}. Using Proposition~\ref{prop:set-soc-convex-1es}, we demonstrate the equivalence between \eqref{eq:multistage} and \eqref{eq:multistage-box-two} as follows:

\begin{proposition}
\label{prop:multi-stage-equivalent-1es}
    Suppose there is only one ESS. Then $(p^{A \vee}, p^{A \wedge})$ is feasible in the Multistage model \eqref{eq:multistage} if and only if there exist $e^{S \downarrow}$ and $e^{S \uparrow}$ such that $(p^{A \vee}, p^{A \wedge}, e^{S \downarrow}, e^{S \uparrow})$ is feasible in the Rectangular model \eqref{eq:multistage-box-two}. Additionally, this holds if and only if there exist $p^{S \vee \downarrow}$, $y^{\vee \downarrow}$, $p^{S \vee \uparrow}$, $y^{\vee \uparrow}$, $p^{S \wedge \downarrow}$, $y^{\wedge \downarrow}$, $p^{S \wedge \uparrow}$, and $y^{\wedge \uparrow}$ such that, together with $(p^{A \vee}, p^{A \wedge}, e^{S \downarrow}, e^{S \uparrow})$, they form a feasible solution to the program \eqref{eq:multistage-equivalent-1es}. Consequently, the optimal values satisfy $\eqref{eq:multistage} = \eqref{eq:multistage-box-two} = \eqref{eq:multistage-equivalent-1es}$.
    \begin{subequations}
    \label{eq:multistage-equivalent-1es}
    \begin{align}
        \textbf{Single-ESS}:\; \max_{\substack{
        p^{A \vee}, p^{A \wedge}, e^{S \downarrow}, e^{S \uparrow}, \\
        p^{S \vee \downarrow}, y^{\vee \downarrow}, p^{S \vee \uparrow}, y^{\vee \uparrow}, \\
        p^{S \wedge \downarrow}, y^{\wedge \downarrow}, p^{S \wedge \uparrow}, y^{\wedge \uparrow}
        }}\; & \phi(p^{A \vee}, p^{A \wedge}) \\
        \text{s.t.}\; & \eqref{eq:pa-order}, \eqref{eq:multistage-box-2}, \eqref{eq:multistage-box-3}, \\
        \label{eq:multistage-equivalent-1es-c}
        & (p_t^{A \vee}, p_t^{S \vee \downarrow}, y_t^{\vee \downarrow}) \in C_t,\, \forall t \in S_T, \\
        \label{eq:multistage-equivalent-1es-d}
        & e_t^{S \downarrow} \leq \kappa^S e_{t - 1}^{S \downarrow} - F_{\eta^{SD}, \eta^{SC}}(p_t^{S \vee \downarrow}) \leq e_t^{S \uparrow},\, \forall t \in S_T, \\
        & (p_t^{A \vee}, p_t^{S \vee \uparrow}, y_t^{\vee \uparrow}) \in C_t,\, \forall t \in S_T, \\
        \label{eq:multistage-equivalent-1es-f}
        & e_t^{S \downarrow} \leq \kappa^S e_{t - 1}^{S \uparrow} - F_{\eta^{SD}, \eta^{SC}}(p_t^{S \vee \uparrow}) \leq e_t^{S \uparrow},\, \forall t \in S_T, \\
        & (p_t^{A \wedge}, p_t^{S \wedge \downarrow}, y_t^{\wedge \downarrow}) \in C_t,\, \forall t \in S_T, \\
        & e_t^{S \downarrow} \leq \kappa^S e_{t - 1}^{S \downarrow} - F_{\eta^{SD}, \eta^{SC}}(p_t^{S \wedge \downarrow}) \leq e_t^{S \uparrow},\, \forall t \in S_T, \\
        & (p_t^{A \wedge}, p_t^{S \wedge \uparrow}, y_t^{\wedge \uparrow}) \in C_t,\, \forall t \in S_T, \\
        \label{eq:multistage-equivalent-1es-j}
        & e_t^{S \downarrow} \leq \kappa^S e_{t - 1}^{S \uparrow} - F_{\eta^{SD}, \eta^{SC}}(p_t^{S \wedge \uparrow}) \leq e_t^{S \uparrow},\, \forall t \in S_T.
    \end{align}
    \end{subequations}
\end{proposition}

The proof of Proposition~\ref{prop:multi-stage-equivalent-1es} is provided in \ref{sec:proof}. We refer to \eqref{eq:multistage-equivalent-1es} as the \emph{Single-ESS} model, which is a single-stage program.

\subsubsection{Disaggregation}

We propose a disaggregation approach for the Rectangular model \eqref{eq:multistage-box-two}. Similar to Algorithm~\ref{alg:enumeration}, we employ a greedy algorithm and minimize the operation cost $c(p_{t_0}^S, y_{t_0})$ for each period $t_0$. Suppose $(p^{A \vee}, p^{A \wedge}, e^{S \downarrow}, e^{S \uparrow})$ is a feasible solution to the Rectangular model \eqref{eq:multistage-box-two}, the optimization problem for disaggregation of period $t_0 \in T$ is formulated as follows:
\begin{subequations}
\label{eq:multistage-box-two-greedy}
\begin{align}
    \min_{e_{t_0}^S, p_{t_0}^S, y_{t_0}}\; & c(p_{t_0}^S, y_{t_0}) \\
    \text{s.t.}\; & \eqref{eq:cons-1}, \eqref{eq:cons-2},\, \text{for}\; t = t_0,  \\
    \label{eq:multistage-box-two-greedy-c}
    & e_{i t_0}^{S \downarrow} \leq e_{i t_0}^S \leq e_{i t_0}^{S \uparrow},\, \forall i \in S_S,
\end{align}
\end{subequations}
where \eqref{eq:multistage-box-two-greedy-c} ensures future feasibility through the SoC range $\times_{i \in S_S} [e_{i t_0}^{S \downarrow}, e_{i t_0}^{S \uparrow}]$. The overall procedure for aggregation and disaggregation is summarized in Algorithm~\ref{alg:rectangular}. Algorithms for the cases of convex ESS constraints or a single ESS can be derived similarly based on \eqref{eq:multistage-box-two-enumeration} and \eqref{eq:multistage-equivalent-1es}.

\begin{algorithm}
    \caption{Aggregation and disaggregation of the Rectangular method}
    \begin{algorithmic}[1]
    \label{alg:rectangular}
        \STATE \textbf{DSO:}
        \STATE \quad Aggregation: Solve \eqref{eq:multistage-box-two} to obtain $p^{A \vee}$, $p^{A \wedge}$, $e^{S \downarrow}$, and $e^{S \uparrow}$.
        \STATE \quad Send $p^{A \vee}$ and $p^{A \wedge}$ to the TSO.
        \FOR{$t_0 = 1$ \TO $T$}
        \STATE \textbf{TSO:}
        \STATE \quad Decide $p_{t_0}^A \in [p_{t_0}^{A \vee}, p_{t_0}^{A \wedge}]$. Send $p_{t_0}^A$ to the DSO.
        \STATE \textbf{DSO:}
        \STATE \quad Disaggregation: Solve \eqref{eq:multistage-box-two-greedy} to obtain $e_{t_0}^S$, $p_{t_0}^S$, and $y_{t_0}$.
        \ENDFOR
    \end{algorithmic}
\end{algorithm}

\subsection{Inner Approximation Based on Envelopes}

This subsection analyzes an envelope-based inner approximation method proposed by \cite{chen2019aggregate} for the two-stage model. We show that this method provides feasible solutions to the Multistage model \eqref{eq:multistage}, though it is more conservative than the Rectangular model \eqref{eq:multistage-box-two}. Additionally, we present an envelope-based greedy disaggregation algorithm.

\subsubsection{Aggregation}

\cite{chen2019aggregate} proposed a single-stage program to inner approximate the two-stage model for time-decoupled power flexibility aggregation, under the assumptions of ideal ESSs, a specific power flow model, and a linear objective function. We adapt it to our general settings (non-ideal ESSs, arbitrary convex power flow models, and general objective functions) while preserving the main ideas, resulting in the following \emph{Envelope} model:
\begin{subequations}
\label{eq:reformulation1-1dim}
\begin{align}
    \label{eq:reformulation1-1dim-obj}
    \textbf{Envelope}:\; \max_{\substack{
        p^{A \vee}, p^{S \vee}, y^\vee, p^{A \wedge}, p^{S \wedge}, y^\wedge, \\
        e^{S \downarrow}, p^{S \downarrow}, e^{S \uparrow}, p^{S \uparrow} 
        }}\; & \phi(p^{A \vee}, p^{A \wedge}) \\
    \label{eq:reformulation1-1dim-soc-initial}
    \text{s.t.}\; & \eqref{eq:pa-order}, e_{i 0}^{S \downarrow} = e_{i 0}^{S \uparrow} = e_{i 0}^S,\, \forall i \in S_S, \\
    \label{eq:reformulation1-1dim-dynamic-down}
    & e_{i t}^{S \uparrow} = \kappa_i^S e_{i (t-1)}^{S \uparrow} - F_{\eta^{SD}_i, \eta^{SC}_i}(p_{i t}^{S \downarrow}),\, \forall i \in S_S,\, \forall t \in S_T,  \\
    \label{eq:reformulation1-1dim-dynamic-up}
    & e_{i t}^{S \downarrow} = \kappa_i^S e_{i (t-1)}^{S \downarrow} - F_{\eta^{SD}_i, \eta^{SC}_i}(p_{i t}^{S \uparrow}),\, \forall i \in S_S,\, \forall t \in S_T, \\
    \label{eq:reformulation1-1dim-convex}
    & (p_t^{A \vee}, p_t^{S \vee}, y_t^\vee), (p_t^{A \wedge}, p_t^{S \wedge}, y_t^\wedge) \in C_t,\, \forall t \in S_T, \\
    \label{eq:reformulation1-1dim-power-bound}
    & p_{i t}^{S \downarrow} \leq p_{i t}^{S \vee} \leq p_{i t}^{S \uparrow},\, p_{i t}^{S \downarrow} \leq p_{i t}^{S \wedge} \leq p_{i t}^{S \uparrow},\, \forall i \in S_S,\, \forall t \in S_T, \\
    \label{eq:reformulation1-1dim-soc-bound}
    & \underline{E}_i^S \leq e_{i t}^{S \downarrow} \leq \overline{E}_i^S,\, \underline{E}_i^S \leq e_{i t}^{S \uparrow} \leq \overline{E}_i^S,\, \forall i \in S_S,\, \forall t \in S_T,
\end{align}
\end{subequations}
where the lower and upper envelopes of the ESS power are represented by $p^{S \downarrow}$ and $p^{S \uparrow}$, respectively. The SoC levels $e^{S \downarrow}$ and $e^{S \uparrow}$ derived from these envelops are required to satisfy the SoC bounds as shown in \eqref{eq:reformulation1-1dim-soc-bound}. The following proposition states that the Envelope model \eqref{eq:reformulation1-1dim} is at least as conservative as the Rectangular model \eqref{eq:multistage-box-two}.

\begin{proposition}
\label{prop:1dim-sufficient}
    Assume $p^{A \vee}$, $p^{S \vee}$, $y^\vee$, $p^{A \wedge}$, $p^{S \wedge}$, $y^\wedge$, $e^{S \downarrow}$, $p^{S \downarrow}$, $e^{S \uparrow}$, and $p^{S \uparrow}$ form a feasible solution to the Envelope model \eqref{eq:reformulation1-1dim}. Then $(p^{A \vee}, p^{A \wedge}, e^{S \downarrow}, e^{S \uparrow})$ is feasible in the Rectangular model \eqref{eq:multistage-box-two}. Consequently, their optimal values satisfy $\eqref{eq:reformulation1-1dim} \leq \eqref{eq:multistage-box-two}$.
\end{proposition}

The proof of Proposition~\ref{prop:1dim-sufficient} is provided in \ref{sec:proof}. Since the Rectangular model \eqref{eq:multistage-box-two} is an inner approximation for the Multistage model \eqref{eq:multistage}, the Envelope model \eqref{eq:reformulation1-1dim} also serves as an inner approximation. Furthermore, \eqref{eq:reformulation1-1dim} is often strictly more conservative than \eqref{eq:multistage-box-two}, as will be illustrated in Section~\ref{sec:case}.

\subsubsection{Disaggregation}

In \eqref{eq:reformulation1-1dim}, the envelopes $p^{S \downarrow}$ and $p^{S \uparrow}$ provide time-decoupled feasible variation region for ESSs. In operation management, larger flexibility intervals are more preferable. Therefore, after solving \eqref{eq:reformulation1-1dim} and obtaining the optimal value $I^*$, we maximize the feasible variation region of ESSs while preserving $I^*$:
\begin{subequations}
\label{eq:maximize-es}
\begin{align}
    \max\; & \sum_{i \in S_S} \phi(p_i^{S \downarrow}, p_i ^{S \uparrow}) \\
    \text{s.t.}\; & \eqref{eq:reformulation1-1dim-soc-initial}\text{--}\eqref{eq:reformulation1-1dim-soc-bound},\, \phi(p^{A \vee}, p^{A \wedge}) \geq I^*.
\end{align}
\end{subequations}

Similar to Algorithm~\ref{alg:enumeration} and Algorithm~\ref{alg:rectangular}, the disaggregation process includes an operation cost minimization program, formulated as follows: For $t_0 \in S_T$,
\begin{subequations}
\label{eq:online}
\begin{align}
    \min_{p_{t_0}^S, y_{t_0}}\; & c(p_{t_0}^S, y_{t_0}) \\
    \text{s.t.}\; & (p_{t_0}^A, p_{t_0}^S, y_{t_0}) \in C_{t_0}, \\
    \label{eq:online-c}
    & p_{i t_0}^{S \downarrow} \leq p_{i t_0}^S \leq p_{i t_0}^{S \uparrow},\, \forall i \in S_S.
\end{align}
\end{subequations}
Constraint \eqref{eq:online-c} ensures that the ESS power levels remain within the envelopes, guaranteeing $e_{i t_0}^{S \downarrow} \leq e_{i t_0}^S \leq e_{i t_0}^{S \uparrow}$ for any $i \in S_S$ and thus ensuring future feasibility. When the cost function $c$ is convex, problem \eqref{eq:online} is convex for each period $t_0 \in S_T$. The procedure for flexibility aggregation and disaggregation in the Envelope method is summarized in Algorithm~\ref{alg:envelope}.

\begin{algorithm}
    \caption{Aggregation and disaggregation of the Envelope method}
    \begin{algorithmic}[1]
    \label{alg:envelope}
        \STATE \textbf{DSO:}
        \STATE \quad Aggregation: Solve \eqref{eq:reformulation1-1dim} to obtain the optimal value $I^*$. Solve \eqref{eq:maximize-es} to obtain $p^{A \vee}$, $p^{A \wedge}$, $p^{S \downarrow}$, and $p^{S \uparrow}$.
        \STATE \quad Send $p^{A \vee}$ and $p^{A \wedge}$ to the TSO.
        \FOR{$t_0 = 1$ \TO $T$}
        \STATE \textbf{TSO:}
        \STATE \quad Decide $p_{t_0}^A \in [p_{t_0}^{A \vee}, p_{t_0}^{A \wedge}]$. Send $p_{t_0}^A$ to the DSO.
        \STATE \textbf{DSO:}
        \STATE \quad Disaggregation: Solve \eqref{eq:online} to obtain $p_{t_0}^S$ and $y_{t_0}$. Calculate $e_{t_0}^S$ by \eqref{eq:cons-2} for $t = t_0$.
        \ENDFOR
    \end{algorithmic}
\end{algorithm}

\subsection{Outer Approximations Based on the Two-Stage Model}

As shown in Proposition~\ref{prop:two-stage}, the Two-Stage model \eqref{eq:two-stage} serves as an outer approximation for the Multistage model \eqref{eq:multistage}. Program \eqref{eq:two-stage} is a two-stage RO problem with DDU. If the convex ESS constraints are employed and the operational constraints are linear, the improved C\&CG algorithm \citep{zeng2022two} can be used to obtain an exact solution for \eqref{eq:two-stage}. However, this approach may be time-consuming occasionally, as two-stage RO is generally NP-hard. To address this, we propose another outer approximation by relaxing constraint \eqref{eq:two-stage-con} as the combination of \eqref{eq:pa-order} and the following constraint:
\begin{align}
    \label{eq:two-stage-con-relax-b}
    & \forall p^A \in \{p^{A \vee}, p^{A \wedge}\},\, \exists e^S, p^S, y,\, \text{s.t.}\; \eqref{eq:cons},\, \forall t \in S_T,
\end{align}
where \eqref{eq:pa-order} ensures that $[p_t^{A \vee}, p_t^{A \wedge}]$ is well-defined, and \eqref{eq:two-stage-con-relax-b} follows from \eqref{eq:two-stage-con} by restricting the uncertainty set to the two elements $p^{A \vee}$ and $p^{A \wedge}$. By reformulating \eqref{eq:two-stage-con-relax-b}, we derive the following single-stage program:
\begin{subequations}
\label{eq:reformulation-outer}
\begin{align}
    \textbf{Outer}:\; \max_{\substack{
        p^{A \vee}, p^{S \vee}, y^\vee, \\
        p^{A \wedge}, p^{S \wedge}, y^\wedge 
        }}\; & \phi(p^{A \vee}, p^{A \wedge}) \\
    \text{s.t.}\; & \eqref{eq:pa-order},\, (p_t^{A \vee}, p_t^{S \vee}, y_t^\vee), (p_t^{A \wedge}, p_t^{S \wedge}, y_t^\wedge) \in C_t,\, \forall t \in S_T, \\
    & e_{i t}^{S \vee} = \kappa_i^S e_{i (t-1)}^{S \vee} - F_{\eta^{SD}_i, \eta^{SC}_i}(p_{i t}^{S \vee}),\, \forall i \in S_S,\, \forall t \in S_T, \\
    & e_{i t}^{S \wedge} = \kappa_i^S e_{i (t-1)}^{S \wedge} - F_{\eta^{SD}_i, \eta^{SC}_i}(p_{i t}^{S \wedge}),\, \forall i \in S_S,\, \forall t \in S_T, \\
    & \underline{E}_i^S \leq e_{i t}^{S \vee} \leq \overline{E}_i^S,\, \underline{E}_i^S \leq e_{i t}^{S \wedge} \leq \overline{E}_i^S,\, \forall i \in S_S,\, \forall t \in S_T.
\end{align}
\end{subequations}
We refer to program \eqref{eq:reformulation-outer} as the \emph{Outer} model, whose relationship to the Two-Stage model \eqref{eq:two-stage} is summarized as follows:

\begin{proposition}
\label{prop:outer-two-stage}
    If $(p^{A \vee}, p^{A \wedge})$ is feasible in the Two-Stage model \eqref{eq:two-stage}, then there exist $p^{S \vee}$, $y^\vee$, $p^{S \wedge}$, and $y^\wedge$ such that $(p^{A \vee}, p^{S \vee}, y^\vee, p^{A \wedge}, p^{S \wedge}, y^\wedge)$ is feasible in the Outer model \eqref{eq:reformulation-outer}. Therefore, the optimal values satisfy $\eqref{eq:two-stage} \leq \eqref{eq:reformulation-outer}$.
\end{proposition}

\subsection{Comparison of Solution Methods}

The characteristics of the aforementioned models are compared in Table~\ref{tab:comparison}. Their relationships are summarized in Theorem~\ref{thm:comparison}, which integrates the results of Proposition~\ref{prop:two-stage}, Corollary~\ref{cor:enumeration}, Corollary~\ref{cor:rectangular}, Proposition~\ref{prop:multi-stage-equivalent-1es}, Proposition~\ref{prop:1dim-sufficient}, and Proposition~\ref{prop:outer-two-stage}.

\begin{table}[!t]
\renewcommand{\arraystretch}{1.2}
\caption{Model Comparison}
\label{tab:comparison}
\centering
\footnotesize
\vspace{0.5em}
\begin{tabular}{cccccc}
\hline
Formulation & Name & Stage & Nonanticipativity & ESS constraints \\ 
\hline
\eqref{eq:two-stage} & Two-Stage & Two-stage & $\times$ & General/convex \\
\eqref{eq:multistage} & Multi-Stage & Multistage & \checkmark & General/convex \\
\eqref{eq:multistage-enumeration} & Enumeration & Single-stage & \checkmark & Convex \\
\eqref{eq:multistage-box-two} & Rectangular & Two-stage & \checkmark & General/convex \\
\eqref{eq:multistage-equivalent-1es} & Single-ESS & Single-stage & \checkmark & General/convex (single ESS) \\
\eqref{eq:reformulation1-1dim} & Envelope & Single-stage & \checkmark & General/convex \\
\eqref{eq:reformulation-outer} & Outer & Single-stage & $\times$ & General/convex \\
\hline
\end{tabular}
\end{table}

\begin{theorem}
\label{thm:comparison}
The optimal values of the models in Table~\ref{tab:comparison} satisfy:
    \begin{enumerate}[(i)]
        \item In general, Envelope $\leq$ Rectangular $\leq$ Multistage $\leq$ Two-Stage $\leq$ Outer.
        \item Under convex ESS constraints, Envelope $\leq$ Rectangular $\leq$ Multistage = Enumeration $\leq$ Two-Stage $\leq$ Outer.
        \item When there is only one ESS, Envelope $\leq$ Single-ESS = Rectangular = Multistage $\leq$ Two-Stage $\leq$ Outer.
        \item When there is no ESS, Envelope = Rectangular = Multistage = Enumeration = Two-Stage = Outer.
    \end{enumerate}
\end{theorem}

\section{Case Studies}
\label{sec:case}

In this section, we first construct simple examples to illustrate the effects of the ESS complementarity constraint and the relationships between different aggregation methods. Subsequently, numerical experiments are conducted using two test cases.

\subsection{Simple Examples}
\label{sec:case-simple}

The following example compares the aggregate flexibility of a non-ideal ESS under different types of constraints.

\begin{example}
\label{eg:complementarity}
    Consider a distribution system with a single node connected to a non-ideal energy storage. The time horizon consists of only one period. The parameters are $\overline{P}_1^{SD} = \overline{P}_1^{SC} = 1$, $\tau = 1$, $\kappa_1^S = 1$, $e_{1, 0}^S = 1$, and $\eta_1^{SD} = \eta_1^{SC} = 0.9$.
    \begin{itemize}
        \item \emph{General ESS constraints}: The operational constraints under general ESS constraints \eqref{eq:ESS} are as follows:
        \begin{subequations}
        \begin{align}
            & e_{1, 1}^S = 1 - \max \{p_{1, 1}^S, 0\} / 0.9 - \min \{p_{1, 1}^S, 0\} \times 0.9, \nonumber \\
            & - 1 \leq p_{1, 1}^S \leq 1,\, 0 \leq e_{1, 1}^S \leq 1,\, p_1^A = - p_{1, 1}^S. \nonumber
        \end{align}
        \end{subequations}
        It is straightforward to verify that the aggregate power flexibility region is $S_{A G} = [-0.9, 0]$.
        \item \emph{Without complementarity constraint}: If complementarity constraint is not considered, the operational constraints are:
        \begin{subequations}
        \label{eq:eg1-no}
        \begin{align}
            & e_{1, 1}^S = 1 - p_{1, 1}^{SD} / 0.9 + p_{1, 1}^{SC} \times 0.9,\, p_1^A = - p_{1, 1}^{SD} + p_{1, 1}^{SC} \\
            & 0 \leq p_{1, 1}^{SD} \leq 1,\, 0 \leq p_{1, 1}^{SC} \leq 1,\, 0 \leq e_{1, 1}^S \leq 1.
        \end{align}
        \end{subequations}
        Then $p_1^A \in S_{A N} = [-0.9, 0.19]$. The upper bound $0.19$ is strictly larger than $0$, achieved by $p_{1, 1}^{SD} = 0.81$ and $p_{1, 1}^{SC} = 1$. This demonstrates that simultaneous charging and discharging of the ESS enlarge the estimated aggregate flexibility.
        \item \emph{Convex ESS constraints}: When convex ESS constraints \eqref{eq:ESS-convex} are applied, the constraint $p_{1, 1}^{SD} + p_{1, 1}^{SC} \leq 1$ is added to \eqref{eq:eg1-no}, and the feasible variation region of $p_1^A$ becomes $S_{A C} = [-0.9, 19/181]$. To achieve $p_1^A = 19/181$, the charging and discharging power values are $p_{1, 1}^{SC} = 100/181$ and $p_{1, 1}^{SD} = 81/181$. These values represent average power levels and can be implemented by charging at power $1$ for time $100/181$ and then discharging at power $1$ for time $81/181$.
    \end{itemize}
    In this example, $S_{A G} \subsetneq S_{A C} \subsetneq S_{A N}$. Thus, neglecting the complementarity constraints of non-ideal ESSs can lead to an overestimation of aggregate flexibility, particularly by overestimating the maximum power demand, as simultaneous charging and discharging of non-ideal ESSs increases the power loss estimation. The convex ESS constraints are able to mitigate this overestimation, and simultaneous charging and discharging can be prevented by treating the power values as averages \citep{shen2020modeling}.
\end{example}

The following two examples show that the Envelope model \eqref{eq:reformulation1-1dim}, the Multistage model \eqref{eq:multistage}, and the Two-Stage model \eqref{eq:two-stage} are not equivalent.

\begin{example}
\label{eg:congestion}
    Consider a distribution system with nodes $1$ and $2$. The time horizon consists of two periods. Node $2$ connects an ideal ESS, a DG, and a load, with the following parameters: $\overline{P}_2^{SD} = \overline{P}_2^{SC} = 1$, $e_{2, 0}^S = 1$, $\kappa_2^S = 1$, $\tau = 1$, $\underline{P}_2^G = (0, 0)$, $\overline{P}_2^G = (1, 0)$, $\underline{P}_2^D = (0, 0)$, and $\overline{P}_2^D = (2, 0)$. Node $1$ has no facilities. The capacity of the transmission line from node $1$ to node $2$ is $1$ and there is no power loss. We consider the aggregate flexibility at node $1$. Thus, $p_t^A = p_{2, t}^D - p_{2, t}^G - p_{2, t}^S$ with $-1 \leq p_t^A \leq 1$ for $t = 1, 2$. The objective function is defined as $\phi(p^{A \vee}, p^{A \wedge}) = (p_1^{A \wedge} - p_1^{A \vee}) + (p_2^{A \wedge} - p_2^{A \vee})$. It can be verified that 
    \begin{align}
        \text{Envelope} = \text{Single-ESS} = \text{Rectangular} = \text{Multistage} = 3 < \text{Two-Stage} = \text{Outer} = 4. \nonumber
    \end{align}
    In fact, any $p^A \in [-1, 1] \times [-1, 1]$ is feasible in the Two-Stage model \eqref{eq:two-stage} by setting $p_{2, 1}^S = (p_2^A + 1)/2$, $p_{2, 1}^D - p_{2, 1}^G = p_1^A + (p_2^A + 1)/2$, $p_{2, 2}^S = - p_2^A$, and $p_{2, 2}^D = p_{2, 2}^G = 0$. However, this solution requires knowledge of $p_2^A$ for the operation in period $1$. Now consider the Multistage model \eqref{eq:multistage}. Without knowing $p_2^A$, the feasible variation range of $p_2^A$ can only have a length of at most $1$, as the ESS is the only flexibility source in period $2$, and its power range is restricted by the SoC bound constraint $0 \leq e_{2, 1}^S - p_{2, 2}^S \leq 1$, regardless of the value of $e_{2, 1}^S$. Therefore, this example demonstrates that the Multistage model and the Two-Stage model are not equivalent in general, even for the case of a single ideal ESS.
\end{example}

\begin{example}
\label{eg:omega}
    Consider a distribution system with a single node, an ideal ESS, and a DG. The time horizon consists of three periods. The parameters are $\overline{P}_1^{SD} = \overline{P}_1^{SC} = 1$, $e_{1, 0}^S = 0$, $\kappa_1^S = 1$, $\tau = 1$, $\underline{P}_1^G = (0, 0, 0)$, and $\overline{P}_1^G = (0, 2, 0)$. The objective function is $\phi(p^{A \vee}, p^{A \wedge}) = 2(p_1^{A \wedge} - p_1^{A \vee}) + (p_2^{A \wedge} - p_2^{A \vee}) + 2(p_3^{A \wedge} - p_3^{A \vee})$. It can be verified that
    \begin{align}
        \text{Envelope} = 4 < \text{Single-ESS} = \text{Rectangular} = \text{Multistage} = 5 < \text{Two-Stage} = \text{Outer} = 6. \nonumber
    \end{align}
    For the Envelope model \eqref{eq:reformulation1-1dim}, the aggregate flexibility in period $1$ (or period $3$) is provided by the ESS, the aggregate flexibility in period $2$ is provided by the DG, and the optimal value is $2 \times 1 + 2 = 4$. For the Multistage model \eqref{eq:multistage}, the aggregate flexibility in period $1$ is provided by the ESS. The flexibility of the DG in period $2$ is divided: Half is allocated to the aggregate flexibility in period $2$, and the other half is used to reset the SoC of the ESS, enabling the ESS to provide flexibility in period $3$. Thus, the optimal value is $2 \times 1 + 1 + 2 \times 1 = 5$ in the Multistage model. Since the weight of period $3$ is higher than that of period $2$, the Multistage model transfers part of the DG's flexibility from period $2$ to period $3$ via the ESS. This example shows that the Envelope model \eqref{eq:reformulation1-1dim}, the Multistage model \eqref{eq:multistage}, and the Two-Stage model \eqref{eq:two-stage} are not equivalent in general, even for the case of a single ideal ESS in a distribution system without congestion.
\end{example}

\subsection{Numerical Experiments}

We evaluate and compare different aggregation and disaggregation methods using numerical experiments. The test systems include modified PJM 5-bus system and IEEE 33-bus system, with parameters provided in \cite{xie2025}. For each system, we examine multiple cases with varying numbers of periods and ESSs. The numerical experiments are based on the DC power flow model, convex ESS constraints \eqref{eq:ESS-convex}, and linear objective functions $\phi(p^{A \vee}, p^{A \wedge}) = \sum_{t \in \mathcal{T}} \omega_t (p_t^{A \wedge} - p_t^{A \vee})$, where $\omega_t$ is a randomly and independently generated coefficient in $[0, 1]$ for each $t \in \mathcal{T}$. Single-stage optimization problems are solved using Gurobi 11.0. Algorithms are coded using Julia. All experiments are conducted on a laptop with an Intel i7-12700H processor and 16 GB RAM.

We first test the aggregation models, including the Envelope, Rectangular, Enumeration, Two-Stage, and Outer models in Table~\ref{tab:comparison}. Additionally, the multistage model with affine policies (referred to as the \emph{Affine} model) is included as a benchmark. The Envelope, Enumeration, and Outer models are solved directly using Gurobi. The Rectangular and Two-Stage models are solved using the C\&CG algorithm with normalization reformulation \citep{zeng2022two}. The Affine model is solved using the improved constraint generation algorithm \citep{su2022multi}. The limit of computation time is set as $10,000$ s. 

The aggregation results across different methods are presented in Table~\ref{tab:aggregation}. We consider $12$ cases with varying number of nodes ($|S_N|$), number of periods ($|S_T|$), and number of ESSs ($|S_S|$), with randomly generated ESSs' parameters in each case. In all cases, the optimal values of the aggregation models satisfy Envelope $\leq$ Rectangular $\leq$ Enumeration $\leq$ Two-Stage $\leq$ Outer, which is consistent with Theorem~\ref{thm:comparison}. Since Enumeration = Multistage under convex ESS constraints, we can see that Envelope $<$ Rectangular and Multistage $<$ Two-Stage often happen. In addition, the Rectangular model improves the flexibility by up to 29.9\% compared to the Envelope model.

In terms of computation efficiency, the Enumeration model cannot be solved for $|S_T| = 24$ due to its large problem size, as the number of variables grows exponentially with $|S_T|$. The Affine model fails to converge within the $10,000$ s time limit for $|S_T| = 24$ and it requires comparable or longer computation time than the Enumeration model for $|S_T| = 8$. In contrast, the Envelope and Outer models are much faster than other methods, as they solve single-stage optimization problems of manageable scales. The Rectangular and Two-Stage models are slower, as they rely on iterative processes to solve two-stage RO problems. The solution time of the Two-Stage model has greater variability, because the number of iterations varies significantly across cases. In comparison, the Rectangular model converges within $3$ iterations in all $12$ cases.

\begin{table}[!t]
\renewcommand{\arraystretch}{1.2}
\caption{Aggregation Results of Different Methods}
\label{tab:aggregation}
\centering
\footnotesize
\vspace{0.5em}
\begin{tabular}{c|c|c|c|c|c|c|c|c}
\hline
\multicolumn{3}{c|}{Settings} & \multicolumn{6}{c}{Optimal value, computation time (s)} \\
\hline
$|S_N|$ & $|S_T|$ & $|S_S|$ & Envelope & Rectangular & Enumeration & Two-Stage & Outer & Affine \\ 
\hline
\multirow{6}{*}{5} & \multirow{3}{*}{8} & 1 & 40.58, 0.075 & 41.89, 0.350 & 41.89, 2.573 & 42.35, 0.150 & 42.35, 0.058 & 41.89, 0.585 \\ 
\cline{3-9}
& & 2 & 40.81, 0.100 & 41.74, 0.478 & 41.74, 2.981 & 41.76, 0.167 & 41.76, 0.058 & 41.74, 1.168 \\ 
\cline{3-9}
& & 5 & 27.43, 0.083 & 32.56, 0.478 & 32.56, 5.610 & 32.56, 0.197 & 32.56, 0.075 & 32.56, 20.94 \\ 
\cline{2-9}
& \multirow{3}{*}{24} & 1 & 98.51, 0.084 & 112.9, 0.108 & -, - & 113.9, 424.1 & 113.9, 0.075 & -, - \\
\cline{3-9}
& & 2 & 95.98, 0.092 & 111.3, 0.961 & -, - & 112.9, 221.8 & 112.9, 0.067 & -, - \\
\cline{3-9}
& & 5 & 97.47, 1.571 & 125.1, 15.31 & -, - & 132.0, 3.101 & 132.0, 0.734 & -, - \\
\hline
\multirow{6}{*}{33} & \multirow{3}{*}{8} & 7 & 28.56, 0.036 & 30.71, 0.192 & 30.71, 6.519 & 31.77, 0.149 & 31.77, 0.017 & 30.71, 7.886 \\ 
\cline{3-9}
& & 16 & 49.78, 0.038 & 56.07, 0.370 & 56.07, 8.635 & 56.07, 0.198 & 56.07, 0.022 & 56.07, 67.50 \\ 
\cline{3-9}
& & 33 & 86.28, 0.646 & 88.72, 7.158 & 88.72, 20.96 & 90.48, 1.567 & 90.48, 0.391 & 88.72, 114.8 \\ 
\cline{2-9}
& \multirow{3}{*}{24} & 7 & 61.97, 0.078 & 79.20, 0.732 & -, - & 83.23, 1.051 & 83.23, 0.146 & -, - \\
\cline{3-9}
& & 16 & 73.34, 0.148 & 95.29, 1.926 & -, - & 95.30, 1.869 & 95.30, 0.098 & -, - \\
\cline{3-9}
& & 33 & 120.7, 0.950 & 148.4, 20.95 & -, - & 150.2, 4.250 & 150.2, 0.414 & -, - \\
\hline
\end{tabular}
\end{table}

\begin{figure}[!t]
\centering
\includegraphics[width=16cm]{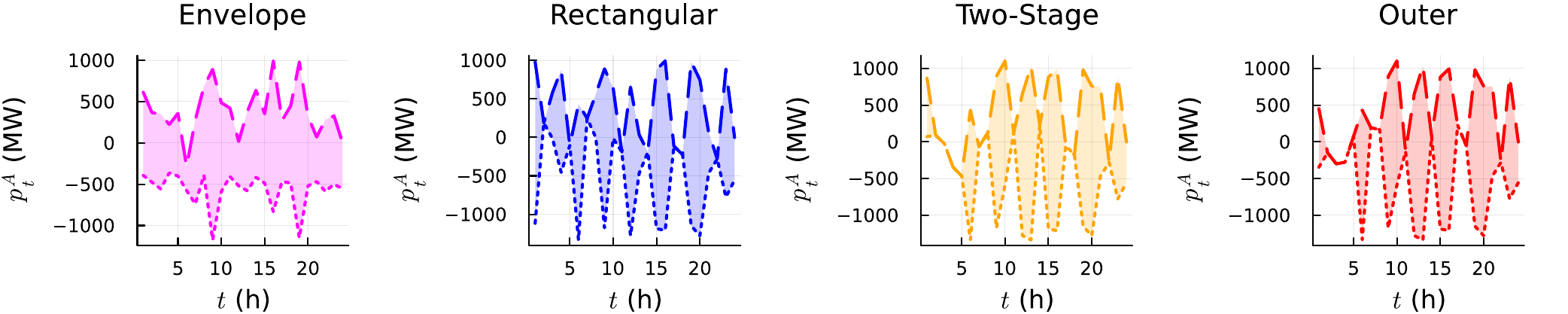}
\caption{Aggregate flexibility in the 5-bus system with $24$ periods and $5$ ESSs. The upper bound $p^{A \wedge}$ and the lower bound $p^{A \vee}$ are depicted by dashed and dotted lines, respectively.}
\label{fig:aggregate_power_case5}
\end{figure}

\begin{figure}[!t]
\centering
\includegraphics[width=16cm]{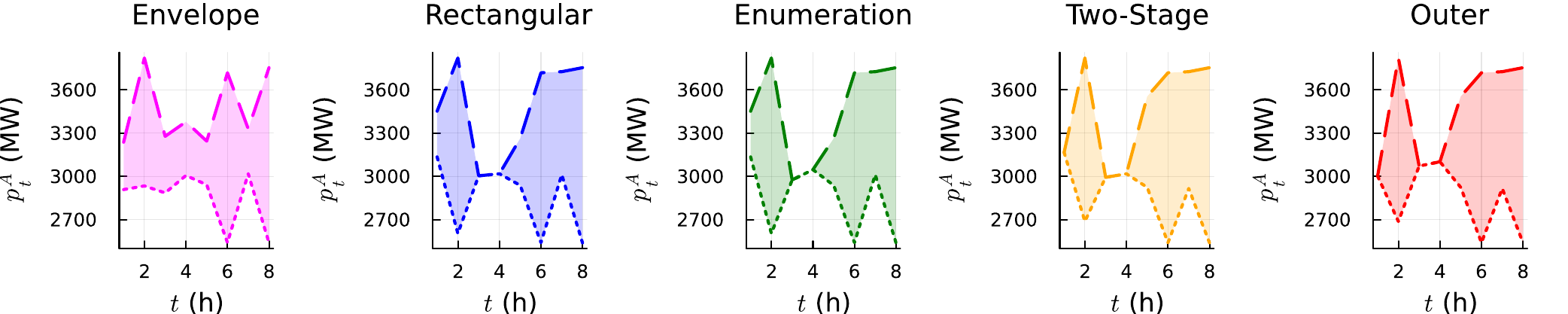}
\caption{Aggregate flexibility in the 33-bus system with $8$ periods and $7$ ESSs. The upper bound $p^{A \wedge}$ and the lower bound $p^{A \vee}$ are depicted by dashed and dotted lines, respectively.}
\label{fig:aggregate_power_case33bw}
\end{figure}

The aggregate power flexibility region, outlined by the curves $p^{A \vee}$ and $p^{A \wedge}$, is illustrated in Figure~\ref{fig:aggregate_power_case5} and Figure~\ref{fig:aggregate_power_case33bw} for two cases. These figures show that different models yield different flexibility aggregation results. In the case settings, DG and load flexibility are present in every period. Thus, $p_t^{A \vee} < p_t^{A \wedge}$ for all $t \in S_T$ under the Envelope model. In contrast, other models allow the transfer of flexibility from periods with a relatively low weight $\omega_t$ to periods with higher weights via ESSs, leading to some periods where $p_t^{A \vee} = p_t^{A \wedge}$.

\begin{figure}[!t]
\centering
\includegraphics[width=16cm]{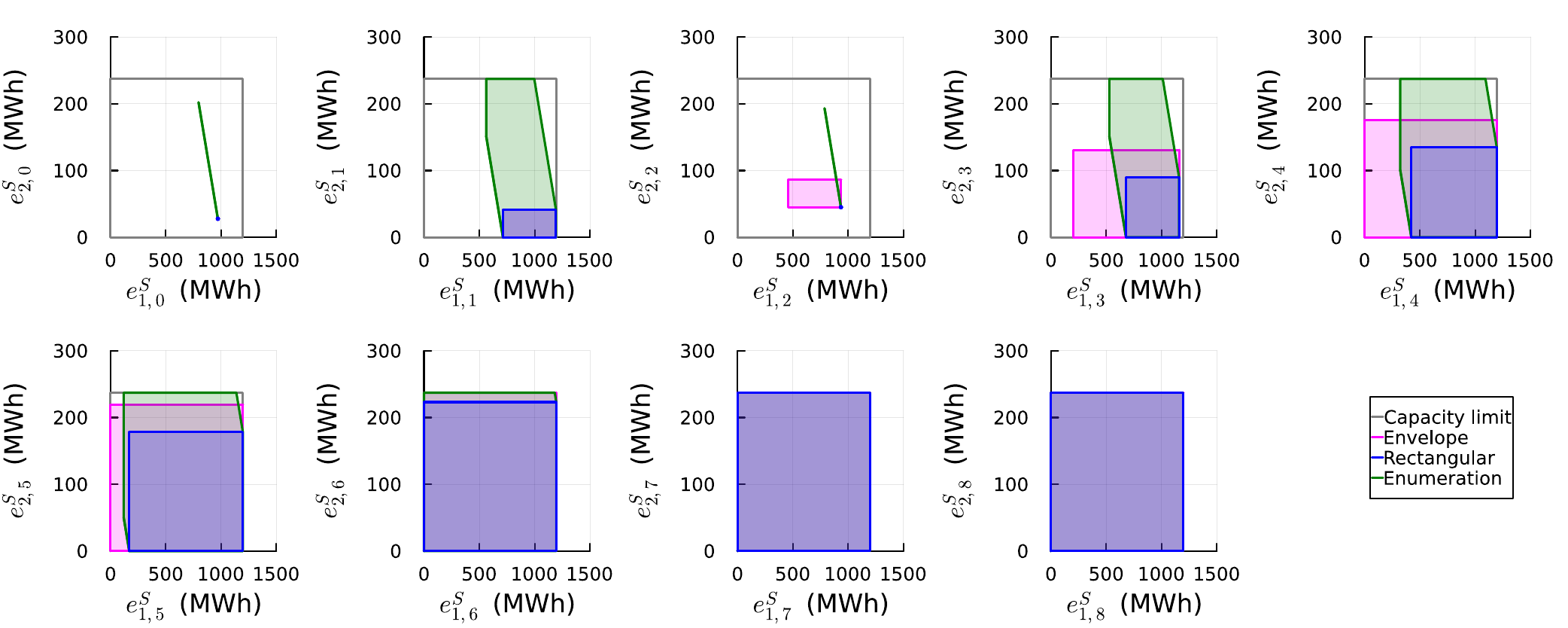}
\caption{SoC ranges in the 5-bus system with $8$ periods and $2$ ESSs.}
\label{fig:SoC}
\end{figure}

To visualize and compare the methods related to SoC ranges, we plot the SoC ranges of the Envelope, Rectangular, and Enumeration models for the 5-bus system with $8$ periods and $2$ ESSs in Figure~\ref{fig:SoC}, where the capacity limits of the SoC ranges are also plotted. The figure verifies that the SoC range in the Envelope and Rectangular models are always box-shaped, whereas those in the Enumeration model can be other polyhedrons. The sloping edges of these polyhedrons have similar slope values, indicating that the SoC levels of the two ESSs can complement each other to some extent in terms of ensuring future flexibility.

Next, we examine the disaggregation results of models incorporating nonanticipativity, including the Envelope, Rectangular, Enumeration, and Affine models. For the Envelope model, we test both the disaggregation method from \cite{chen2019aggregate} and Algorithm~\ref{alg:envelope}. For each case, $100$ aggregate power trajectories $p^A$ are randomly generated within the intersection of the aggregate power flexibility regions derived by these models, and we compute the average operation cost and computation time. The disaggregation results are summarized in Table~\ref{tab:disaggregation}. Algorithm~\ref{alg:envelope}, the Rectangular model, and the Enumeration model optimize the DSO's operation strategy to minimize operation costs, resulting in lower average operation costs compared to the other two methods that do not consider the operation costs in the decision process. In particular, Algorithm 3 lowers the average operation costs by up to 40.3\% compared to the method in \cite{chen2019aggregate}. Optimization is not involved in the disaggregation method in \cite{chen2019aggregate} and the Affine model, so they are computationally fast. In contrast, Algorithm~\ref{alg:envelope} and the Rectangular model are slower due to the need to solve single-stage optimization problems, but they are still fast enough for real-time operation management. The aggregation and disaggregation of the Enumeration model share a similar computational burden, and this method is the most time-consuming in Table~\ref{tab:disaggregation}.

\begin{table}[!t]
\renewcommand{\arraystretch}{1.2}
\caption{Disaggregation Results of Different Methods}
\label{tab:disaggregation}
\centering
\footnotesize
\vspace{0.5em}
\begin{tabular}{c|c|c|c|c|c|c|c}
\hline
\multicolumn{3}{c|}{Settings} & \multicolumn{5}{c}{Average operation cost (\$$10^5$), average computation time (s)} \\
\hline
$|S_N|$ & $|S_T|$ & $|S_S|$ & Envelope & Envelope & Rectangular & Enumeration & Affine \\ 
& & & \citep{chen2019aggregate} & (Algorithm~\ref{alg:envelope}) & & & \\
\hline
\multirow{6}{*}{5} & \multirow{3}{*}{8} & 1 & 1644, 0.000 & 1409, 0.116 & 1360, 0.116 & 1355, 20.54 & 1912, 0.006 \\ 
\cline{3-8}
& & 2 & 1590, 0.000 & 1314, 0.117 & 1290, 0.117 & 1279, 23.84 & 1855, 0.000 \\ 
\cline{3-8}
& & 5 & 1635, 0.000 & 1357, 0.116 & 1256, 0.117 & 1248, 33.89 & 1830, 0.000 \\ 
\cline{2-8}
& \multirow{3}{*}{24} & 1 & 4933, 0.000 & 4198, 0.349 & 4356, 0.348 & -, - & -, - \\
\cline{3-8}
& & 2 & 4835, 0.000 & 4116, 0.346 & 4276, 0.353 & -, - & -, - \\
\cline{3-8}
& & 5 & 4884, 0.001 & 4137, 0.358 & 4150, 0.357 & -, - & -, - \\
\hline
\multirow{6}{*}{33} & \multirow{3}{*}{8} & 7 & 9.531, 0.000 & 5.691, 0.007 & 5.702, 0.004 & 5.701, 49.87 & 9.777, 0.007 \\ 
\cline{3-8}
& & 16 & 9.455, 0.000 & 5.694, 0.010 & 5.733, 0.008 & 5.684, 72.26 & 9.230, 0.000 \\ 
\cline{3-8}
& & 33 & 8.971, 0.001 & 5.905, 0.020 & 5.661, 0.015 & 5.661, 104.3 & 10.01, 0.000 \\ 
\cline{2-8}
& \multirow{3}{*}{24} & 7 & 30.18, 0.000 & 20.02, 0.018 & 20.67, 0.014 & -, - & -, - \\
\cline{3-8}
& & 16 & 30.33, 0.000 & 20.09, 0.023 & 20.00, 0.020 & -, - & -, - \\
\cline{3-8}
& & 33 & 30.49, 0.002 & 20.03, 0.037 & 19.91, 0.034 & -, - & -, - \\
\hline
\end{tabular}
\end{table}

\section{Conclusion}
\label{sec:conclusion}

This paper proposes a multistage RO model for the DSO's time-decoupled power flexibility aggregation, where the TSO's decision for the aggregate power trajectory is treated as sequentially revealed DDU. Unlike the existing two-stage model where the disaggregation process assumes full knowledge of the aggregate trajectory, our model has the advantage of considering the TSO's sequential decision-making process. We propose various methods to solve the multistage model exactly or approximately under different assumptions of ESS constraints. Corresponding disaggregation methods based on greedy algorithms are also developed. The application scope, nonanticipativity, conservativeness, and computational burden of these methods are analyzed and compared both theoretically and through case studies. 

The main conclusions are as follows: 1) The Multistage model is generally more conservative than the Two-Stage model due to the incorporation of nonanticipativity constraints. 2) When convex ESS constraints are applied, the Enumeration model can solve the Multistage model exactly; however, it becomes impractical for cases with a large number of periods, as the number of variables and constraints grows exponentially. 3) For cases with general ESS constraints, the Envelope and Rectangular models serve as inner approximations for the Multistage model. The Envelope model is faster due to its single-stage form, while the Rectangular model offers more accurate approximations but is more computationally demanding, as it requires iterative solution algorithms for its two-stage structure. Nonetheless, the iteration number for the Rectangular model is typically below $3$ in the numerical experiments. The Two-Stage and Outer models act as outer approximations for the Multistage model. The Outer model is single-stage and computationally efficient, whereas the Two-Stage model is more accurate but requires a longer solution time. 4) Incorporating the greedy algorithm and considering operation cost minimization in disaggregation effectively reduces the average operation cost, as demonstrated by the numerical experiments. 5) The complementarity constraint preventing simultaneous ESS charging and discharging cannot be neglected in power flexibility aggregation for non-ideal ESSs. Applying convex ESS constraints helps maintain convexity, facilitating the solution process.

Based on these findings, we recommend the following strategies for time-decoupled power flexibility aggregation: 1) For ideal ESSs, the complementarity constraint can be neglected, while convex ESS constraints are recommended for non-ideal ESSs to maintain convexity. 2) When the number of periods is small, the Enumeration model is suitable for exact solutions to the Multistage model. For larger time horizons, the Envelope and Outer models can be used to quickly obtain inner and outer approximations, respectively. If additional computational time is available, the Rectangular and Two-Stage models provide more accurate approximations. 3) To reduce operation costs, the greedy disaggregation algorithms can be employed, combined with the aggregation methods incorporating nonanticipativity.

Several directions for future work may further enhance the applicability. First, applying the model predict control technique and incorporating the uncertainty of DGs and loads in the power flexibility aggregation could improve the adaptability of the models to real-world variability. Second, exploring the use of an aggregate cost function would enable more economically efficient operation strategies at the transmission system level. Third, leveraging parallel computing could improve scalability, allowing the models to handle larger systems with greater efficiency.

\bibliographystyle{apalike}
\bibliography{mybib}

\newpage

\appendix

\section{Proofs}
\label{sec:proof}

\begin{proof}[Proof of Lemma~\ref{lemma:ESS-our-common}]
    Suppose $(e_{i (t - 1)}^S, e_{i t}^S, p_{i t}^S)$ satisfies \eqref{eq:ESS-dynamic} and \eqref{eq:ESS-p}. Let $p_{i t}^{SD} = \max \{p_{i t}^S, 0\}$ and $p_{i t}^{SC} = - \min \{p_{i t}^S, 0\}$. Because $- \overline{P}_i^{SC} < 0 < \overline{P}_i^{SD}$, we have $0 \leq p_{i t}^{SD} \leq \overline{P}_i^{SD}$ and $0 \leq p_{i t}^{SC} \leq \overline{P}_i^{SC}$. In addition, $p_{i t}^{SD} - p_{i t}^{SC} = \max \{p_{i t}^S, 0\} + \min \{p_{i t}^S, 0\} = p_{i t}^S$ and $p_{i t}^{SD} p_{i t}^{SC} = \max \{p_{i t}^S, 0\} \min \{p_{i t}^S, 0\} = 0$. The SoC dynamic constraint also holds. Therefore, $(e_{i (t - 1)}^S, e_{i t}^S, p_{i t}^S)$ is in \eqref{eq:ESS-common}.

    Suppose $(e_{i (t - 1)}^S, e_{i t}^S, p_{i t}^S)$ is in \eqref{eq:ESS-common}, then there exist $p_{i t}^{SD}$ and $p_{i t}^{SC}$ that satisfy the constraints of \eqref{eq:ESS-common}. Because at most one of $p_{i t}^{SD}$ and $p_{i t}^{SC}$ can be nonzero, $p_{i t}^S = p_{i t}^{SD}$ if $p_{i t}^{SD} > 0$; $p_{i t}^S = - p_{i t}^{SC}$ if $p_{i t}^{SC} > 0$; $p_{i t}^S = 0$ if $p_{i t}^{SD} = p_{i t}^{SC} = 0$. Therefore, $p_{i t}^{SD} = \max \{p_{i t}^S, 0\}$ and $p_{i t}^{SC} = - \min \{p_{i t}^S, 0\}$. Then it is easy to see that $(e_{i (t - 1)}^S, e_{i t}^S, p_{i t}^S)$ satisfies \eqref{eq:ESS-dynamic} and \eqref{eq:ESS-p}.
\end{proof}

\begin{proof}[Proof of Lemma~\ref{lemma:ESS-ideal}]
    By Lemma~\ref{lemma:ESS-our-common}, we only need to prove that \eqref{eq:ESS-common} and \eqref{eq:ESS-relaxed} are equal as sets. When $\eta_i^{SD} = \eta_i^{SC} = 1$, relaxing the complementarity constraint in \eqref{eq:ESS-common} gives \eqref{eq:ESS-relaxed}. Thus, \eqref{eq:ESS-common} is a subset of \eqref{eq:ESS-relaxed}. Suppose $(e_{i (t - 1)}^S, e_{i t}^S, p_{i t}^S)$ is in \eqref{eq:ESS-relaxed}. Then there exist $p_{i t}^{SD}$ and $p_{i t}^{SC}$ satisfy the constraints of \eqref{eq:ESS-relaxed}, so $- \overline{P}_i^{SC} \leq p_{i t}^S \leq \overline{P}_i^{SD}$. Let $\tilde{p}_{i t}^{SD} = \max \{p_{i t}^S, 0\}$ and $\tilde{p}_{i t}^{SC} = - \min \{p_{i t}^S, 0\}$. Then $p_{i t}^S = \tilde{p}_{i t}^{SD} - \tilde{p}_{i t}^{SC}$, $0 \leq \tilde{p}_{i t}^{SD} \leq \overline{P}_i^{SD}$, $0 \leq \tilde{p}_{i t}^{SC} \leq \overline{P}_i^{SC}$, $\tilde{p}_{i t}^{SD} \tilde{p}_{i t}^{SC} = 0$ and $e_{i t}^S = \kappa_i^S e_{i (t-1)}^S - p_{i t}^{SD} \tau + p_{i t}^{SC} \tau = \kappa_i^S e_{i (t-1)}^S - p_{i t}^S \tau = \kappa_i^S e_{i (t-1)}^S - \tilde{p}_{i t}^{SD} \tau + \tilde{p}_{i t}^{SC} \tau$. Therefore, $(e_{i (t - 1)}^S, e_{i t}^S, p_{i t}^S, \tilde{p}_{i t}^{SD}, \tilde{p}_{i t}^{SC})$ satisfy the constraints of \eqref{eq:ESS-common}. Therefore, \eqref{eq:ESS-relaxed} equals \eqref{eq:ESS-common}.
\end{proof}

\begin{proof}[Proof of Proposition~\ref{prop:set-soc-reformulate}]
    We reformulate $\mathcal{E}_{t_0}^S(p^{A \vee}, p^{A \wedge})$ and \eqref{eq:set-soc-reformulate} using decision rules. For $t_0 < t \in S_T$, let $p_{(t_0, t]}^A = (p_{t_0 + 1}^A, p_{t_0 + 2}^A, \dots, p_t^A)$ denote the vector of revealed aggregate power after period $t_0$, and $p_{> t}^A = (p_{t + 1}^A, \dots, p_T^A)$ denote the future aggregate power. Let $e_t^S(e_{t_0}^S, p_{(t_0, t]}^A)$ represent the decision rule for $e_t^S = (e_{i t}^S; i \in S_S)$, which depends on the SoC $e_{t_0}^S$ at the end of period $t_0$ and the revealed aggregate power $p_{(t_0, t]}^A$, but is independent of $p_{> t}^A$. Operation decision rules $p_t^{SD}(e_{t_0}^S, p_{(t_0, t]}^A)$, $p_t^{SC}(e_{t_0}^S, p_{(t_0, t]}^A)$, and $y_t(e_{t_0}^S, p_{(t_0, t]}^A)$ are similar. With these notations, $\mathcal{E}_{t_0}^S(p^{A \vee}, p^{A \wedge})$ and \eqref{eq:set-soc-reformulate} can be rewritten as follows:
    \begin{align}
        & \mathcal{E}_{t_0}^S(p^{A \vee}, p^{A \wedge}) = \left\{ e_{t_0}^S \,\middle|\,  
        \begin{aligned}
            & \underline{E}_i^S \leq e_{i t_0}^S \leq \overline{E}_i^S,\, \forall i \in S_S, \\
            & \forall p_{> t_0}^A \in \times_{t_0 < t \in S_T} [p_t^{A \vee}, p_t^{A \wedge}], \\
            & \exists e_t^S(e_{t_0}^S, p_{(t_0, t]}^A), p_t^{SD} (e_{t_0}^S, p_{(t_0, t]}^A), p_t^{SC} (e_{t_0}^S, p_{(t_0, t]}^A), \\
            & y_t (e_{t_0}^S, p_{(t_0, t]}^A); t_0 < t \in S_T,\, \text{s.t.}\; \eqref{eq:cons-convex},\, \forall t > t_0
        \end{aligned}
        \right\}, \nonumber \\
        & \eqref{eq:set-soc-reformulate} = \left\{ e_{t_0}^S ~\middle|~  
        \begin{aligned}
            & \underline{E}_i^S \leq e_{i t_0}^S \leq \overline{E}_i^S,\, \forall i \in S_S, \\
            & \forall p_{> t_0}^A \in \times_{t_0 < t \in S_T} \{p_t^{A \vee}, p_t^{A \wedge}\}, \\
            & \exists e_t^S(e_{t_0}^S, p_{(t_0, t]}^A), p_t^{SD} (e_{t_0}^S, p_{(t_0, t]}^A), p_t^{SC} (e_{t_0}^S, p_{(t_0, t]}^A), \\
            & y_t (e_{t_0}^S, p_{(t_0, t]}^A); t_0 < t \in S_T,\, \text{s.t.}\; \eqref{eq:cons-convex},\, \forall t > t_0
        \end{aligned}
        \right\}. \nonumber
    \end{align}

    Then $\eqref{eq:set-soc-reformulate} \supseteq \mathcal{E}_{t_0}^S(p^{A \vee}, p^{A \wedge})$ because $\times_{t_0 < t \in S_T} \{p_t^{A \vee}, p_t^{A \wedge}\}$ is a subset of $\times_{t_0 < t \in S_T} [p_t^{A \vee}, p_t^{A \wedge}]$. To prove $\eqref{eq:set-soc-reformulate} \subseteq \mathcal{E}_{t_0}^S(p^{A \vee}, p^{A \wedge})$, assume $e_{t_0}^S \in \eqref{eq:set-soc-reformulate}$. For any $p_{> t_0}^A \in \times_{t_0 < t \in S_T} [p_t^{A \vee}, p_t^{A \wedge}]$, there exists $\lambda \in \times_{t_0 < t \in S_T} [0, 1]$ such that $p_t^A = \lambda_t p_t^{A \vee} + (1 - \lambda_t) p_t^{A \wedge}$ for any $t_0 < t \in S_T$. Let
    \begin{subequations}
    \begin{align}
        e_{t_0 + 1}^S(e_{t_0}^S, p_{(t_0, t_0 + 1]}^A) =\, & \lambda_{t_0 + 1} e_{t_0 + 1}^S(e_{t_0}^S, p_{t_0 + 1}^{A \vee}) + (1 - \lambda_{t_0 + 1}) e_{t_0 + 1}^S(e_{t_0}^S, p_{t_0 + 1}^{A \wedge}), \nonumber \\
        e_{t_0 + 2}^S(e_{t_0}^S, p_{(t_0, t_0 + 2]}^A) =\, & \lambda_{t_0 + 1} (\lambda_{t_0 + 2} e_{t_0 + 2}^S (e_{t_0}^S, p_{t_0 + 1}^{A \vee}, p_{t_0 + 2}^{A \vee}) \nonumber \\
        & + (1 - \lambda_{t_0 + 2}) e_{t_0 + 2}^S (e_{t_0}^S, p_{t_0 + 1}^{A \vee}, p_{t_0 + 2}^{A \wedge})) \nonumber \\
        & + (1 - \lambda_{t_0 + 1}) (\lambda_{t_0 + 2} e_{t_0 + 2}^S (e_{t_0}^S, p_{t_0 + 1}^{A \wedge}, p_{t_0 + 2}^{A \vee}) \nonumber \\
        & + (1 - \lambda_{t_0 + 2}) e_{t_0 + 2}^S (e_{t_0}^S, p_{t_0 + 1}^{A \wedge}, p_{t_0 + 2}^{A \wedge})), \nonumber \\
        e_{t_0 + 3}^S (e_{t_0}^S, p_{(t_0, t_0 + 3]}^A) =\, & \lambda_{t_0 + 1} (\lambda_{t_0 + 2} ( \lambda_{t_0 + 3} e_{t_0 + 3}^S (e_{t_0}^S, p_{t_0 + 1}^{A \vee}, p_{t_0 + 2}^{A \vee}, p_{t_0 + 3}^{A \vee}) \nonumber \\
        & + (1 - \lambda_{t_0 + 3}) e_{t + 3}^S (e_{t_0}^S, p_{t_0 + 1}^{A \vee}, p_{t_0 + 2}^{A \vee}, p_{t_0 + 3}^{A \wedge})) \nonumber \\
        & + (1 - \lambda_{t_0 + 2}) ( \lambda_{t_0 + 3} e_{t_0 + 3}^S (e_{t_0}^S, p_{t_0 + 1}^{A \vee}, p_{t_0 + 2}^{A \wedge}, p_{t_0 + 3}^{A \vee}) \nonumber \\
        & + (1 - \lambda_{t_0 + 3}) e_{t_0 + 3}^S (e_{t_0}^S, p_{t_0 + 1}^{A \vee}, p_{t_0 + 2}^{A \wedge}, p_{t_0 + 3}^{A \wedge}))) \nonumber \\
        & + (1 - \lambda_{t_0 + 1}) (\lambda_{t_0 + 2} ( \lambda_{t_0 + 3} e_{t_0 + 3}^S (e_{t_0}^S, p_{t_0 + 1}^{A \wedge}, p_{t_0 + 2}^{A \vee}, p_{t_0 + 3}^{A \vee}) \nonumber \\
        & + (1 - \lambda_{t_0 + 3}) e_{t_0 + 3}^S (e_{t_0}^S, p_{t_0 + 1}^{A \wedge}, p_{t_0 + 2}^{A \vee}, p_{t_0 + 3}^{A \wedge})) \nonumber \\
        & + (1 - \lambda_{t_0 + 2}) ( \lambda_{t_0 + 3} e_{t + 3}^S (e_{t_0}^S, p_{t_0 + 1}^{A \wedge}, p_{t_0 + 2}^{A \wedge}, p_{t_0 + 3}^{A \vee}) \nonumber \\
        & + (1 - \lambda_{t_0 + 3}) e_{t_0 + 3}^S (e_{t_0}^S, p_{t_0 + 1}^{A \wedge}, p_{t_0 + 2}^{A \wedge}, p_{t_0 + 3}^{A \wedge}))), \nonumber \\
        \dots, & \nonumber
    \end{align}
    \end{subequations}
    where the right-hand side is well-defined because $e_{t_0}^S \in \eqref{eq:set-soc-reformulate}$. The terms $p_t^{SD}(e_{t_0}^S, p_{(t_0, t]}^A)$, $p_t^{SC}(e_{t_0}^S, p_{(t_0, t]}^A)$, and $y_t(e_{t_0}^S, p_{(t_0, t]}^A)$ can be defined similarly for $t_0 < t \in S_T$. To demonstrate that the constraints in $\mathcal{E}_{t_0}^S(p^{A \vee}, p^{A \wedge})$ are satisfied, it suffices to verify \eqref{eq:cons-convex} for $t_0 < t \in S_T$: Constraints \eqref{eq:cons-convex-1} and \eqref{eq:cons-convex-2} follow from the convexity of $C_t$ and the convex ESS constraints, while \eqref{eq:cons-convex-3} follows from the convexity of the set $\{ e_t^S ~|~ \underline{E}_i^S \leq e_{i t}^S \leq \overline{E}_i^S, \forall i \in S_S \}$.
\end{proof}

\begin{proof}[Proof of Proposition~\ref{prop:set-soc-equivalent}]
    When $(p^{A \vee}, p^{A \wedge})$ is feasible in \eqref{eq:multistage}, let $\mathcal{E}_t^S = \mathcal{E}_t^S(p^{A \vee}, p^{A \wedge})$ as defined in Definition~\ref{def:SoC-range}. By Lemma~\ref{lemma:set-soc}, these sets are nonempty and closed, and they satisfy $\mathcal{E}_t^S \subseteq \{ e_{i t}^S \,|\, \underline{E}_i^S \leq e_{i t}^S \leq \overline{E}_i^S, \forall i \in S_S \}$ for all $t \in S_T \cup \{0\}$, with $e_0^S \in \mathcal{E}_0^S(p^{A \vee}, p^{A \wedge})$. The expression for $\mathcal{E}_{t - 1}^S(p^{A \vee}, p^{A \wedge})$ in Lemma~\ref{lemma:soc-range} yields \eqref{eq:set-soc-equivalent}.

    Conversely, assume there exist nonempty sets $\mathcal{E}_t^S \subseteq \{ e_t^S \,|\, \underline{E}_i^S \leq e_{i t}^S \leq \overline{E}_i^S, \forall i \in S_S \}$ for all $t \in S_T \cup \{0\}$, satisfying $e_0^S \in \mathcal{E}_0^S$ and \eqref{eq:set-soc-equivalent} for all $t \in S_T$. Since $e_0^S \in \mathcal{E}_0^S$, for any $p_1^A \in [p_1^{A \vee}, p_1^{A \wedge}]$, there exist $e_1^S \in \mathcal{E}_1^S$, $p_1^S$, and $y_1$ such that \eqref{eq:cons} holds for $t = 1$. Similarly, since $e_1^S \in \mathcal{E}_1^S$, for any $p_2^A \in [p_2^{A \vee}, p_2^{A \wedge}]$, there exist $e_2^S \in \mathcal{E}_2^S$, $p_2^S$, and $y_2$ such that \eqref{eq:cons} holds for $t = 2$. Repeat this process up to $t = T$. Then $(p^{A \vee}, p^{A \wedge})$ is feasible in \eqref{eq:multistage}.
\end{proof}

\begin{proof}[Proof of Proposition~\ref{prop:set-soc-convex}]
    To demonstrate the convexity of $\mathcal{E}_t^S(p^{A \vee}, p^{A \wedge})$ for all $t \in S_T \cup \{0\}$, we employ mathematical induction. For $t = T$, the set
    \begin{align}
        \mathcal{E}_T^S(p^{A \vee}, p^{A \wedge}) = \left\{ e_T^S \,\middle|\, \underline{E}_i^S \leq e_{i T}^S \leq \overline{E}_i^S, \forall i \in S_S \right\}, \nonumber
    \end{align}
    is convex. In what follows, we assume $\mathcal{E}_t^S(p^{A \vee}, p^{A \wedge})$ is convex for $t \in S_T$ and prove the convexity of $\mathcal{E}_{t - 1}^S(p^{A \vee}, p^{A \wedge})$.
    
    Assume $\tilde{e}_{t - 1}^S, \hat{e}_{t - 1}^S \in \mathcal{E}_{t - 1}^S(p^{A \vee}, p^{A \wedge})$ and $\lambda \in [0, 1]$. Let $e_{t - 1}^S = \lambda \tilde{e}_{t - 1}^S + (1 - \lambda) \hat{e}_{t - 1}^S$. We need to prove that $e_{t - 1}^S \in \mathcal{E}_{t - 1}^S(p^{A \vee}, p^{A \wedge})$. 
    
    Using the expression for $\mathcal{E}_{t - 1}^S(p^{A \vee}, p^{A \wedge})$ in Lemma~\ref{lemma:soc-range}, for $p_t^A = p_t^{A \vee}$, there exist $\tilde{e}_t^{S \vee}$, $\tilde{p}_t^{SD \vee}$, $\tilde{p}_t^{SC \vee}$, $\tilde{y}_t^\vee$, $\hat{e}_t^{S \vee}$, $\hat{p}_t^{SD \vee}$, $\hat{p}_t^{SC \vee}$, and $\hat{y}_t^\vee$ such that
    \begin{subequations}
    \label{eq:set-soc-convex-vee}
    \begin{align}
        & \tilde{e}_t^{S \vee}, \hat{e}_t^{S \vee} \in \mathcal{E}_t^S(p^{A \vee}, p^{A \wedge}), \\
        \label{eq:set-soc-convex-vee-b}
        & (p_t^{A \vee}, \tilde{p}_t^{SD \vee}, \tilde{p}_t^{SC \vee}, \tilde{y}_t^\vee), (p_t^{A \vee}, \hat{p}_t^{SD \vee}, \hat{p}_t^{SC \vee}, \hat{y}_t^\vee) \in C_t, \\
        & \tilde{e}_{i t}^{S \vee} = \kappa_i^S \tilde{e}_{i (t - 1)}^S - \tilde{p}_{i t}^{SD \vee} \tau / \eta_i^{SD} + \tilde{p}_{i t}^{SC \vee} \tau \eta_i^{SC},\, \forall i \in S_S, \\
        & \hat{e}_{i t}^{S \vee} = \kappa_i^S \hat{e}_{i (t - 1)}^S - \hat{p}_{i t}^{SD \vee} \tau / \eta_i^{SD} + \hat{p}_{i t}^{SC \vee} \tau \eta_i^{SC},\, \forall i \in S_S.
    \end{align}
    \end{subequations}
    Similarly, there exist $\tilde{e}_t^{S \wedge}$, $\tilde{p}_t^{SD \wedge}$, $\tilde{p}_t^{SC \wedge}$, $\tilde{y}_t^\wedge$, $\hat{e}_t^{S \wedge}$, $\hat{p}_t^{SD \wedge}$, $\hat{p}_t^{SC \wedge}$, and $\hat{y}_t^\wedge$ such that
    \begin{subequations}
    \begin{align}
        & \tilde{e}_t^{S \wedge}, \hat{e}_t^{S \wedge} \in \mathcal{E}_t^S(p^{A \wedge}, p^{A \wedge}), \\
        \label{eq:set-soc-convex-wedge-b}
        & (p_t^{A \wedge}, \tilde{p}_t^{SD \wedge}, \tilde{p}_t^{SC \wedge}, \tilde{y}_t^\wedge), (p_t^{A \wedge}, \hat{p}_t^{SD \wedge}, \hat{p}_t^{SC \wedge}, \hat{y}_t^\wedge) \in C_t, \\
        & \tilde{e}_{i t}^{S \wedge} = \kappa_i^S \tilde{e}_{i (t - 1)}^S - \tilde{p}_{i t}^{SD \wedge} \tau / \eta_i^{SD} + \tilde{p}_{i t}^{SC \wedge} \tau \eta_i^{SC},\, \forall i \in S_S, \\
        & \hat{e}_{i t}^{S \wedge} = \kappa_i^S \hat{e}_{i (t - 1)}^S - \hat{p}_{i t}^{SD \wedge} \tau / \eta_i^{SD} + \hat{p}_{i t}^{SC \wedge} \tau \eta_i^{SC},\, \forall i \in S_S.
    \end{align}
    \end{subequations}
    
    For any $\mu \in [0, 1]$ and $p_t^A = \mu p_t^{A \vee} + (1 - \mu) p_t^{A \wedge}$, define
    \begin{subequations}
    \begin{align}
        & e_t^S = \mu (\lambda \tilde{e}_t^{S \vee} + (1 - \lambda) \hat{e}_t^{S \vee}) + (1 - \mu) (\lambda \tilde{e}_t^{S \wedge} + (1 - \lambda) \hat{e}_t^{S \wedge}), \nonumber \\
        & p_t^{SD} = \mu (\lambda \tilde{p}_t^{SD \vee} + (1 - \lambda) \hat{p}_t^{SD \vee}) + (1 - \mu) (\lambda \tilde{p}_t^{SD \wedge} + (1 - \lambda) \hat{p}_t^{SD \wedge}), \nonumber \\
        & p_t^{SC} = \mu (\lambda \tilde{p}_t^{SC \vee} + (1 - \lambda) \hat{p}_t^{SC \vee}) + (1 - \mu) (\lambda \tilde{p}_t^{SC \wedge} + (1 - \lambda) \hat{p}_t^{SC \wedge}), \nonumber \\
        & y_t = \mu (\lambda \tilde{y}_t^\vee + (1 - \lambda) \hat{y}_t^\vee) + (1 - \mu) (\lambda \tilde{y}_t^\wedge + (1 - \lambda) \hat{y}_t^\wedge). \nonumber
    \end{align}
    \end{subequations}
    Then $e_t^S \in \mathcal{E}_t^S(p^{A \vee}, p^{A \wedge})$ follows from the convexity of $\mathcal{E}_t^S(p^{A \vee}, p^{A \wedge})$ and
    \begin{align}
        e_{i t}^S = \kappa_i^S e_{i (t - 1)}^S - p_{i t}^{SD} \tau / \eta_i^{SD} + p_{i t}^{SC} \tau \eta_i^{SC},\, \forall i \in S_S. \nonumber
    \end{align}
    According to \eqref{eq:set-soc-convex-vee-b} and the convexity of $C_t$, we have $(p_t^{A \vee}, \lambda \tilde{p}_t^{SD \vee} + (1 - \lambda) \hat{p}_t^{SD \vee}, \lambda \tilde{p}_t^{SC \vee} + (1 - \lambda) \hat{p}_t^{SC \vee}, \lambda \tilde{y}_t^\vee + (1 - \lambda) \hat{y}_t^\vee) \in C_t$. Similarly, by \eqref{eq:set-soc-convex-wedge-b}, we have $(p_t^{A \wedge}, \lambda \tilde{p}_t^{SD \wedge} + (1 - \lambda) \hat{p}_t^{SD \wedge}, \lambda \tilde{p}_t^{SC \wedge} + (1 - \lambda) \hat{p}_t^{SC \wedge}, \lambda \tilde{y}_t^\wedge + (1 - \lambda) \hat{y}_t^\wedge) \in C_t$. Combining these results, we obtain $(p_t^A, p_t^{SD}, p_t^{SC}, y_t) \in C_t$. Since $\tilde{e}_{t - 1}^S, \hat{e}_{t - 1}^S \in \mathcal{E}_{t - 1}^S(p^{A \vee}, p^{A \wedge})$ and $\lambda \in [0, 1]$ are arbitrary, it follows that $\mathcal{E}_{t - 1}^S(p^{A \vee}, p^{A \wedge})$ is convex. Consequently, by mathematical induction, $\mathcal{E}_t^S(p^{A \vee}, p^{A \wedge})$ is convex for all $t \in S_T \cup \{0\}$.
\end{proof}

\begin{proof}[Proof of Proposition~\ref{prop:multistage-box-uncertainty}]
    Constraint \eqref{eq:set-soc-equivalent-extreme} can be viewed as a relaxation of \eqref{eq:multistage-box-4}, where $\forall p_t^A \in [p_t^{A \vee}, p_t^{A \wedge}]$ is replaced with $\forall p_t^A \in \{p_t^{A \vee}, p_t^{A \wedge}\}$. Consequently, \eqref{eq:multistage-box-4} implies \eqref{eq:set-soc-equivalent-extreme}. 
    
    Conversely, assume that constraint \eqref{eq:set-soc-equivalent-extreme} holds. Fix an arbitrary $t \in S_T$ and $e_{t - 1}^S$ such that
    \begin{align}
        e_{i (t - 1)}^{S \downarrow} \leq e_{i (t - 1)}^S \leq e_{i (t - 1)}^{S \uparrow}, \forall i \in S_S \nonumber.
    \end{align}
    Then, there exist $e_t^{S \vee}$, $p_t^{S \vee}$, $y_t^\vee$, $e_t^{S \wedge}$, $p_t^{S \wedge}$, and $y_t^\wedge$ satisfying the requirements in \eqref{eq:set-soc-equivalent-extreme}. For any $p_t^A \in [p_t^{A \vee}, p_t^{A \wedge}]$, there exists $\lambda \in [0, 1]$ such that $p_t^A = \lambda p_t^{A \vee} + (1 - \lambda) p_t^{A \wedge}$. Define
    \begin{subequations}
    \begin{align}
        & p_t^S = \lambda p_t^{S \vee} + (1 - \lambda) p_t^{S \wedge}, \nonumber \\
        & y_t = \lambda y_t^\vee + (1 - \lambda) y_t^\wedge, \nonumber \\
        & e_{i t}^S = \kappa_i^S e_{i (t - 1)}^S - F_{\eta_i^{SD}, \eta_i^{SC}} (p_{i t}^S),\, \forall i \in S_S. \nonumber
    \end{align}
    \end{subequations}
    By the convexity of $C_t$, we have $(p_t^A, p_t^S, y_t) \in C_t$. It remains to show that $e_{i t}^{S \downarrow} \leq e_{i t} \leq e_{i t}^{S \uparrow}$ for all $i \in S_S$. Using the function $F$ in Definition~\ref{def:es-power} and Lemma~\ref{lemma:es-power}, it follows that
    \begin{align}
        e_{i t}^{S \downarrow} & \leq \min \{ e_{i t}^{S \vee}, e_{i t}^{S \wedge} \} \nonumber \\
        & = \kappa_i^S e_{i (t - 1)}^S - F_{\eta_i^{SD}, \eta_i^{SC}} (\max \{ p_{i t}^{S \vee}, p_{i t}^{S \wedge} \}) \nonumber \\
        & \leq \kappa_i^S e_{i (t - 1)}^S - F_{\eta_i^{SD}, \eta_i^{SC}} (p_{i t}^S) = e_{i t}^S,\, \forall i \in S_S. \nonumber
    \end{align}
    Similarly, $e_{i t}^S \leq e_{i t}^{S \uparrow}$ for all $i \in S_S$. Thus, \eqref{eq:set-soc-equivalent-extreme} and \eqref{eq:multistage-box-4} are equivalent.
\end{proof}

\begin{proof}[Proof of Proposition~\ref{prop:set-soc-convex-1es}]
    We use mathematical induction to demonstrate the convexity of the set $\mathcal{E}_t^S(p^{A \vee}, p^{A \wedge})$ for all $t \in S_T \cup \{0\}$. Similar to the proof of Proposition~\ref{prop:set-soc-convex}, $\mathcal{E}_T^S(p^{A \vee}, p^{A \wedge})$ is convex. Next, we assume $\mathcal{E}_t^S(p^{A \vee}, p^{A \wedge})$ is convex for $t \in S_T$ and prove the convexity of $\mathcal{E}_{t - 1}^S(p^{A \vee}, p^{A \wedge})$.

    Since there is only one ESS, $\mathcal{E}_t^S(p^{A \vee}, p^{A \wedge}) \subseteq \mathbb{R}$ and we omit the subscript for the ESS. Furthermore, $\mathcal{E}_t^S(p^{A \vee}, p^{A \wedge})$ is convex, closed, and bounded. Therefore, we can denote it by $\mathcal{E}_t^S(p^{A \vee}, p^{A \wedge}) = [e_t^{S \downarrow}, e_t^{S \uparrow}]$. Let
    \begin{subequations}
    \begin{align}
        & p_t^{S \vee \downarrow} = \min \{p_t^S \,|\, \exists y_t, \text{s.t.}\; (p_t^{A \vee}, p_t^S, y_t) \in C_t\}, \nonumber \\
        & p_t^{S \wedge \downarrow} = \min \{p_t^S \,|\, \exists y_t, \text{s.t.}\; (p_t^{A \wedge}, p_t^S, y_t) \in C_t\}, \nonumber \\
        & p_t^{S \vee \uparrow} = \max \{p_t^S \,|\, \exists y_t, \text{s.t.}\; (p_t^{A \vee}, p_t^S, y_t) \in C_t\}, \nonumber \\
        & p_t^{S \wedge \uparrow} = \max \{p_t^S \,|\, \exists y_t, \text{s.t.}\; (p_t^{A \wedge}, p_t^S, y_t) \in C_t\}. \nonumber
    \end{align}
    \end{subequations}
    Assume $\check{e}_{t - 1}^S, \hat{e}_{t - 1}^S \in \mathcal{E}_{t - 1}^S(p^{A \vee}, p^{A \wedge})$ with $\check{e}_{t - 1}^S \leq \hat{e}_{t - 1}^S$. By the definition of $\mathcal{E}_{t - 1}^S(p^{A \vee}, p^{A \wedge})$, there exist $\check{p}_t^{S \vee}$ and $\check{y}_t^\vee$ such that $(\check{p}_t^{A \vee}, \check{p}_t^{S \vee}, \check{y}_t^\vee) \in C_t$ and $e_t^{S \downarrow} \leq \kappa^S \check{e}_{t - 1}^S - F_{\eta^{SD}, \eta^{SC}} (\check{p}_t^{S \vee}) \leq e_t^{S \uparrow}$. From the definition of $p_t^{S \vee \uparrow}$, it follows that $p_t^{S \vee \downarrow} \leq \check{p}_t^{S \vee} \leq p_t^{S \vee \uparrow}$. Since the function $F$ defined in Definition~\ref{def:es-power} is strictly increasing, we have $\kappa^S \check{e}_{t - 1}^S - F_{\eta^{SD}, \eta^{SC}}(p_t^{S \vee \downarrow}) \geq \kappa^S \check{e}_{t - 1}^S - F_{\eta^{SD}, \eta^{SC}}(\check{p}_t^{S \vee}) \geq e_t^{S \downarrow}$ and $\kappa^S \check{e}_{t - 1}^S - F_{\eta^{SD}, \eta^{SC}}(p_t^{S \vee \uparrow}) \leq \kappa^S \check{e}_{t - 1}^S - F_{\eta^{SD}, \eta^{SC}}(\check{p}_t^{S \vee}) \leq e_t^{S \uparrow}$. Similarly, we have:
    \begin{subequations}
    \begin{align}
        & \kappa^S \check{e}_{t - 1}^S - F_{\eta^{SD}, \eta^{SC}}(p_t^{S \vee \downarrow}) \geq e_t^{S \downarrow}, \nonumber \\
        & \kappa^S \check{e}_{t - 1}^S - F_{\eta^{SD}, \eta^{SC}}(p_t^{S \vee \uparrow}) \leq e_t^{S \uparrow}, \nonumber \\
        & \kappa^S \check{e}_{t - 1}^S - F_{\eta^{SD}, \eta^{SC}}(p_t^{S \wedge \downarrow}) \geq e_t^{S \downarrow}, \nonumber \\
        & \kappa^S \check{e}_{t - 1}^S - F_{\eta^{SD}, \eta^{SC}}(p_t^{S \wedge \uparrow}) \leq e_t^{S \uparrow}, \nonumber \\
        & \kappa^S \hat{e}_{t - 1}^S - F_{\eta^{SD}, \eta^{SC}}(p_t^{S \vee \downarrow}) \geq e_t^{S \downarrow}, \nonumber \\
        & \kappa^S \hat{e}_{t - 1}^S - F_{\eta^{SD}, \eta^{SC}}(p_t^{S \vee \uparrow}) \leq e_t^{S \uparrow}, \nonumber \\
        & \kappa^S \hat{e}_{t - 1}^S - F_{\eta^{SD}, \eta^{SC}}(p_t^{S \wedge \downarrow}) \geq e_t^{S \downarrow}, \nonumber \\
        & \kappa^S \hat{e}_{t - 1}^S - F_{\eta^{SD}, \eta^{SC}}(p_t^{S \wedge \uparrow}) \leq e_t^{S \uparrow}. \nonumber
    \end{align}
    \end{subequations}
    
    Suppose $e_{t - 1}^S$ satisfies $\check{e}_{t - 1}^S \leq e_{t - 1}^S \leq \hat{e}_{t - 1}^S$. Then,
    \begin{subequations}
    \label{eq:set-soc-convex-1es}
    \begin{align}
        \label{eq:set-soc-convex-1es-1}
        & \kappa^S e_{t - 1}^S - F_{\eta^{SD}, \eta^{SC}}(p_t^{S \vee \uparrow}) \leq \kappa^S \hat{e}_{t - 1}^S - F_{\eta^{SD}, \eta^{SC}}(p_t^{S \vee \uparrow}) \leq e_t^{S \uparrow}, \\
        \label{eq:set-soc-convex-1es-2}
        & \kappa^S e_{t - 1}^S- F_{\eta^{SD}, \eta^{SC}}(p_t^{S \wedge \uparrow}) \leq \kappa^S \hat{e}_{t - 1}^S - F_{\eta^{SD}, \eta^{SC}}(p_t^{S \wedge \uparrow}) \leq e_t^{S \uparrow}, \\
        \label{eq:set-soc-convex-1es-3}
        & \kappa^S e_{t - 1}^S - F_{\eta^{SD}, \eta^{SC}}(p_t^{S \vee \downarrow}) \geq \kappa^S \check{e}_{t - 1}^S - F_{\eta^{SD}, \eta^{SC}}(p_t^{S \vee \downarrow}) \geq e_t^{S \downarrow}, \\
        \label{eq:set-soc-convex-1es-4}
        & \kappa^S e_{t - 1}^S - F_{\eta^{SD}, \eta^{SC}}(p_t^{S \wedge \downarrow}) \geq \kappa^S \check{e}_{t - 1}^S - F_{\eta^{SD}, \eta^{SC}}(p_t^{S \wedge \downarrow}) \geq e_t^{S \downarrow}. 
    \end{align}
    \end{subequations}
    Note that $e_t^{S \downarrow} \leq e_t^{S \uparrow}$ and $\kappa^S e_{t - 1}^S - F_{\eta^{SD}, \eta^{SC}}(p_t^{S \vee \uparrow}) \leq \kappa^S e_{t - 1}^S - F_{\eta^{SD}, \eta^{SC}}(p_t^{S \vee \downarrow})$. Combining these with \eqref{eq:set-soc-convex-1es-1} and \eqref{eq:set-soc-convex-1es-3}, it can be verified that $[e_t^{S \downarrow}, e_t^{S \uparrow}] \cap [\kappa^S e_{t - 1}^S - F_{\eta^{SD}, \eta^{SC}}(p_t^{S \vee \uparrow}), \kappa^S e_{t - 1}^S - F_{\eta^{SD}, \eta^{SC}}(p_t^{S \vee \downarrow})] \neq \varnothing$. Thus, there exists $p_t^{S \vee} \in [p_t^{S \vee \downarrow}, p_t^{S \vee \uparrow}]$ such that $\kappa^S e_{t - 1}^S - F_{\eta^{SD}, \eta^{SC}}(p_t^{S \vee}) \in [e_t^{S \downarrow}, e_t^{S \uparrow}]$. By the definitions of $p_t^{S \vee \downarrow}$ and $p_t^{S \vee \uparrow}$ and the convexity of $C_t$, there exists $y_t^\vee$ such that $(p_t^{A \vee}, p_t^{S \vee}, y_t^\vee) \in C_t$. The case for $p_t^{A \wedge}$ is analogous. Consequently, there exist $e_t^{S \vee}$, $p_t^{S \vee}$, $y_t^\vee$, $e_t^{S \wedge}$, $p_t^{S \wedge}$, and $y_t^\wedge$ such that
    \begin{subequations}
    \begin{align}
        & (p_t^{A \vee}, p_t^{S \vee}, y_t^\vee) \in C_t,\, (p_t^{A \wedge}, p_t^{S \wedge}, y_t^\wedge) \in C_t, \nonumber \\
        & e_t^{S \vee} = \kappa^S e_{t - 1}^S - F_{\eta^{SD}, \eta^{SC}}(p_t^{S \vee}) \in [e_t^{S \downarrow}, e_t^{S \uparrow}], \nonumber \\
        & e_t^{S \wedge} = \kappa^S e_{t - 1}^S - F_{\eta^{SD}, \eta^{SC}}(p_t^{S \wedge}) \in [e_t^{S \downarrow}, e_t^{S \uparrow}]. \nonumber
    \end{align}
    \end{subequations}   

    For any $\mu \in [0, 1]$ and $p_t^A = \mu p_t^{A \vee} + (1 - \mu) p_t^{A \wedge}$, define
    \begin{align}
        p_t^S = \mu p_t^{S \vee} + (1 - \mu) p_t^{S \wedge},\, y_t = \mu y_t^\vee + (1 - \mu) y_t^\wedge. \nonumber
    \end{align}
    Then $(p_t^A, p_t^S, y_t) \in C_t$. Define $e_t^S$ as
    \begin{align}
        e_t^S = \kappa^S e_{t - 1}^S - F_{\eta^{SD}, \eta^{SC}}(p_t^S). \nonumber
    \end{align}
    Thus, by the monotonicity of function $F$,
    \begin{align}
        \kappa^S e_{t - 1}^S - F_{\eta^{SD}, \eta^{SC}}(p_t^S) \in [\min \{e_t^{S \vee}, e_t^{S \wedge}\}, \max \{e_t^{S \vee}, e_t^{S \wedge}\}] \subseteq [e_t^{S \downarrow}, e_t^{S \uparrow}]. \nonumber
    \end{align}
    Therefore, $e_{t - 1}^S \in \mathcal{E}_{t - 1}^S(p^{A \vee}, p^{A \wedge})$. Thus, $\mathcal{E}_{t - 1}^S(p^{A \vee}, p^{A \wedge})$ is convex. By mathematical induction, the proof is complete.
\end{proof}

\begin{proof}[Proof of Proposition~\ref{prop:multi-stage-equivalent-1es}]
    Corollary~\ref{cor:multi-stage-equivalent} establishes the equivalence between \eqref{eq:multistage} and \eqref{eq:multistage-equivalent}. By Proposition~\ref{prop:set-soc-convex-1es}, we can require that $\mathcal{E}_t^S$ is convex in \eqref{eq:multistage-equivalent}. Since there is only one ESS and its SOC is bounded, $\mathcal{E}_t^S$ is bounded, closed, and convex in $\mathbb{R}$. Consequently, $\mathcal{E}_t^S$ must be an interval and can be expressed as $[e_t^{S \downarrow}, e_t^{S \uparrow}]$. Thus, \eqref{eq:multistage} and \eqref{eq:multistage-box-two} are equivalent when there is a single ESS.

    To demonstrate the equivalence between \eqref{eq:multistage-box-two} and \eqref{eq:multistage-equivalent-1es}, observe that
    \begin{subequations}
    \begin{align}
        \label{eq:multi-stage-equivalent-1es-proof1}
        \eqref{eq:set-soc-equivalent-extreme} \iff &
        \left\{
        \begin{aligned}
            & \forall t \in S_T,\, \forall e_{t - 1}^S \in \left[ e_{t - 1}^{S \downarrow}, e_{t - 1}^{S \uparrow} \right], \\
            & \exists e_t^{S \vee}, p_t^{S \vee}, y_t^\vee, e_t^{S \wedge}, p_t^{S \wedge}, y_t^\wedge,\, \text{s.t.}\; \eqref{eq:set-soc-equivalent-extreme-1}\text{--}\eqref{eq:set-soc-equivalent-extreme-4}
        \end{aligned}
        \right. \\
        \implies & \left\{
        \begin{aligned}
            & \forall t \in S_T,\, \forall e_{t - 1}^S \in \left\{ e_{t - 1}^{S \downarrow}, e_{t - 1}^{S \uparrow} \right\}, \\
            & \exists e_t^{S \vee}, p_t^{S \vee}, y_t^\vee, e_t^{S \wedge}, p_t^{S \wedge}, y_t^\wedge,\, \text{s.t.}\; \eqref{eq:set-soc-equivalent-extreme-1}\text{--}\eqref{eq:set-soc-equivalent-extreme-4}
        \end{aligned}
        \right. \\
        \iff & \eqref{eq:multistage-equivalent-1es-c}\text{--}\eqref{eq:multistage-equivalent-1es-j}. \nonumber
    \end{align}
    \end{subequations}
    Thus, it suffices to prove that constraints \eqref{eq:multistage-equivalent-1es-c}--\eqref{eq:multistage-equivalent-1es-j} imply the right-hand side of \eqref{eq:multi-stage-equivalent-1es-proof1}. Assume \eqref{eq:multistage-equivalent-1es-c}--\eqref{eq:multistage-equivalent-1es-j} hold. For any $t \in S_T$ and $e_{t - 1}^{S \downarrow} \leq e_{t - 1}^S \leq e_{t - 1}^{S \uparrow}$, there exists $\lambda \in [0, 1]$ such that $e_{t - 1}^S = \lambda e_{t - 1}^{S \downarrow} + (1 - \lambda) e_{t - 1}^{S \uparrow}$. Define
    \begin{align}
        e_t^{S \vee} & = \lambda (\kappa^S e_{t - 1}^{S \downarrow} - F_{\eta^{SD}, \eta^{SC}}(p_t^{S \vee \downarrow})) + (1 - \lambda)(\kappa^S e_{t - 1}^{S \uparrow} - F_{\eta^{SD}, \eta^{SC}}(p_t^{S \vee \uparrow})) \nonumber \\
        & = \kappa^S e_{t - 1}^S - (\lambda F_{\eta^{SD}, \eta^{SC}}(p_t^{S \vee \downarrow}) + (1 - \lambda) F_{\eta^{SD}, \eta^{SC}}(p_t^{S \vee \uparrow})). \nonumber
    \end{align}
    Then $e_t^{S \downarrow} \leq e_t^{S \vee} \leq e_t^{S \uparrow}$ by \eqref{eq:multistage-equivalent-1es-d} and \eqref{eq:multistage-equivalent-1es-f}. Define
    \begin{align}
        p_t^{S \vee} = F_{\eta^{SD}, \eta^{SC}}^{-1}(\lambda F_{\eta^{SD}, \eta^{SC}}(p_t^{S \vee \downarrow}) + (1 - \lambda) F_{\eta^{SD}, \eta^{SC}}(p_t^{S \vee \uparrow})). \nonumber
    \end{align}
    Then $e_t^{S \vee} = \kappa^S e_{t - 1}^S - F_{\eta^{SD}, \eta^{SC}}(p_t^{S \vee})$. Furthermore, by the monotonicity of $F$, there exists $\mu \in [0, 1]$ such that $p_t^{S \vee} = \mu p_t^{S \vee \downarrow} + (1 - \mu) p_t^{S \vee \uparrow}$. Define
    \begin{align}
        y_t^\vee = \mu p_t^{S \vee \downarrow} + (1 - \mu) p_t^{S \vee \uparrow}. \nonumber
    \end{align}
    Then $e_t^{S \vee}$, $p_t^{S \vee}$, and $y_t^\vee$ satisfy the constraints in \eqref{eq:set-soc-equivalent-extreme-1}--\eqref{eq:set-soc-equivalent-extreme-4}. Similarly, $e_t^{S \wedge}$, $p_t^{S \wedge}$, and $y_t^\wedge$ can be constructed. Consequently, \eqref{eq:multistage-equivalent-1es-c}--\eqref{eq:multistage-equivalent-1es-j} imply the right-hand side of \eqref{eq:multi-stage-equivalent-1es-proof1}, and thus \eqref{eq:multistage}, \eqref{eq:multistage-box-two}, and \eqref{eq:multistage-equivalent-1es} are equivalent.
\end{proof}

\begin{proof}[Proof of Proposition~\ref{prop:1dim-sufficient}]
    Constraints \eqref{eq:pa-order} and \eqref{eq:multistage-box-2} in \eqref{eq:multistage-box-two} follow directly from \eqref{eq:reformulation1-1dim-soc-initial}. By \eqref{eq:reformulation1-1dim-soc-bound}, it suffices to show $e_{i t}^{S \downarrow} \leq e_{i t}^{S \uparrow}$ for all $i \in S_S$ and $t \in S_T$ to prove \eqref{eq:multistage-box-3}. Using \eqref{eq:reformulation1-1dim-soc-initial}--\eqref{eq:reformulation1-1dim-dynamic-up} and the monotonicity of $F$, we have for all $i \in S_S$ and $t \in S_T$,
    \begin{align}
        e_{i t}^{S \uparrow} & = (\kappa_i^S)^t e_{i 0}^S - \sum_{s = 1}^t (\kappa_i^S)^{t - s} F_{\eta_i^{SD}, \eta_i^{SC}}(p_{i s}^{S \downarrow}) \nonumber \\
        & \geq (\kappa_i^S)^t e_{i 0}^S - \sum_{s = 1}^t (\kappa_i^S)^{t - s} F_{\eta_i^{SD}, \eta_i^{SC}}(p_{i s}^{S \uparrow}) = e_{i t}^{S \downarrow}. \nonumber
    \end{align}
    Thus, \eqref{eq:multistage-box-3} holds. Next, we prove \eqref{eq:set-soc-equivalent-extreme}. Fix $t \in S_T$. For any $e_{t - 1}^S \in \times_{i \in S_S} [e_{i (t - 1)}^{S \downarrow}, e_{i (t - 1)}^{S \uparrow}]$, define $e_{i t}^{S \vee}$ and $e_{i t}^{S \wedge}$ according to \eqref{eq:set-soc-equivalent-extreme-3} and \eqref{eq:set-soc-equivalent-extreme-4}. Constraint \eqref{eq:set-soc-equivalent-extreme-1} is identical to \eqref{eq:reformulation1-1dim-convex}. Therefore, it suffices to verify \eqref{eq:set-soc-equivalent-extreme-2}. Using \eqref{eq:reformulation1-1dim-dynamic-up} and \eqref{eq:reformulation1-1dim-power-bound}, we obtain
    \begin{align}
        e_{i t}^{S \vee} & = \kappa_i^S e_{i (t - 1)}^S - F_{\eta_i^{SD}, \eta_i^{SC}} (p_{i t}^{S \vee}) \nonumber \\
        & \geq \kappa_i^S e_{i (t - 1)}^{S \downarrow} - F_{\eta_i^{SD}, \eta_i^{SC}} (p_{i t}^{S \uparrow}) = e_{i t}^{S \downarrow},\, \forall i \in S_S. \nonumber
    \end{align}
    Similarly, $e_{i t}^{S \vee} \leq e_{i t}^{S \uparrow}$, $e_{i t}^{S \wedge} \geq e_{i t}^{S \downarrow}$, and $e_{i t}^{S \wedge} \leq e_{i t}^{S \uparrow}$ can be proven for all $i \in S_S$. Therefore, \eqref{eq:set-soc-equivalent-extreme} holds and $(p^{A \vee}, p^{A \wedge}, e^{S \downarrow}, e^{S \uparrow})$ is feasible in \eqref{eq:multistage-box-two}.
\end{proof}

\end{document}